\documentclass[11pt]{article}
\usepackage{mathrsfs}
\usepackage{epic,eepic,epsf,epsfig}
\usepackage{amsfonts,srcltx,mathrsfs}
\usepackage{multirow}
\usepackage{amsmath}
\usepackage{amssymb}
\usepackage{amsbsy}
\usepackage{graphicx}
\usepackage{amsfonts}
\usepackage{color}

\newtheorem{lem}{Lemma}[section]%
\newtheorem{theorem}[lem]{Theorem}%
\newtheorem{exam}[lem]{Example}%
\newtheorem{prop}[lem]{Proposition}%

\def\a{\alpha} \def\b{\beta} \def\g{\gamma} \def\d{\delta}

 \def\O{\Omega} \def\G{\Gamma}

\def\di{\bigm|} \def\lg{\langle} \def\rg{\rangle}
\def\nd{\mathrel{\bigm|\kern-.7em/}}

\def\f{\noindent}

\def\PSL{\hbox{\rm PSL}}
 \def\PGL{\hbox{\rm PGL}}
\def\P\GammaL{\hbox{\rm P\GammaL}}\def\Mult{\hbox{\rm Mult}}

\def\Aut{\hbox{\rm Aut}}
\def\Dip{\hbox{\rm Dip}}

\def\soc{\hbox{\rm soc}}

\def\Cay{\hbox{\rm Cay}}
\def\BiCay{\hbox{\rm BiCay}}

\def\mod{\hbox{\rm mod }}

\def\Dih{\hbox{\rm Dih}}

\def\mz{{\mathbb Z}}

\def\Dip{\hbox{\rm Dip}}

\begin{document}
\title{Arc-transitive cyclic and dihedral covers of pentavalent symmetric graphs of order twice a prime}

\footnotetext[1]{ Corresponding author. E-mails: dwyang$@$bjtu.edu.cn,
yqfeng$@$bjtu.edu.cn, jxzhou$@$bjtu.edu.cn}

\author{ \\ Da-Wei Yang, Yan-Quan Feng\footnotemark, Jin-Xin Zhou\\
{\small\em Mathematics, Beijing Jiaotong University, Beijing,
100044, P.R. China}\\}

\date{}
 \maketitle

\begin{abstract}

A regular cover of a connected graph is called {\em cyclic} or {\em
dihedral} if its transformation group
is cyclic or dihedral respectively, and {\em arc-transitive}
(or {\em symmetric}) if the fibre-preserving
automorphism subgroup acts arc-transitively on the regular cover. In this paper, we give a classification of arc-transitive cyclic and dihedral covers of a connected pentavalent symmetric graph of order twice a prime.
All those covers are explicitly constructed as Cayley graphs on some groups, and their full automorphism groups are determined.

\bigskip
\f {\bf Keywords:} Symmetric graph, Cayley graph, Bi-Cayley graph, Regular cover.

\end{abstract}

\section{Introduction}

All groups and graphs considered in this paper are finite,
and all graphs are simple, connected and undirected, unless otherwise stated.
Let $G$ be a permutation group on a set $\O$ and let $\a\in \O$.
Denote by $G_{\a}$ the stabilizer of
$\a$ in $G$, that is, the subgroup of $G$ fixing the point $\a$.
We say that $G$ is {\it semiregular} on
$\O$ if $G_{\a}=1$ for every $\a\in\O$, and {\it regular}
if $G$ is transitive and semiregular. Denote by $\mz_n$, $\mz_n^*$,
$D_n$, $A_n$ and $S_n$ the cyclic group of order $n$,
the multiplicative group of units of $\mz_n$, the dihedral group of order $2n$,
the alternating and symmetric group of degree $n$, respectively.
For two groups $M$ and $N$, we use $MN$, $M.N$, $M\rtimes N$ and $M\times N$ to denote
the product of $M$ and $N$, an extension of $M$ by $N$, a split extension of $M$ by $N$ and
the direct product of $M$ and $N$, respectively.
For a subgroup $H$ of a group $G$, $C_G(H)$ means the centralizer of $H$ in $G$ and $N_G(H)$ means the normalizer of $H$ in $G$.

For a graph $\G$, we denote its vertex set,
edge set and full automorphism group by $V(\G)$, $E(\G)$ and $\Aut(\G)$,
respectively. An {\it $s$-arc} in $\G$ is an ordered
$(s+1)$-tuple $(v_0,v_1,\ldots,v_s)$ of vertices of $\G$ such that
$\{v_{i-1},v_i\}\in E(\G)$ for $1\leq i\leq s$, and $v_{i-1}\neq
v_{i+1}$ for $1\leq i<s$. A $1$-arc is just an {\it arc}. A graph $\G$ is
{\it $(G,s)$-arc-transitive} if $G\leq \Aut(\G)$ acts transitively
on the set of $s$-arcs of $\G$, and {\it $(G, s)$-transitive} if $\G$
is $(G,s)$-arc-transitive but not $(G,s+1)$-arc-transitive.
A graph $\G$ is said to be {\it $s$-arc-transitive} or {\it $s$-transitive} if it is
$(\Aut(\G), s)$-arc-transitive or $(\Aut(\G), s)$-transitive, respectively.
In particular, 0-arc-transitive means {\it vertex-transitive}, and 1-arc-transitive means {\it arc-transitive} or {\it symmetric}.
A graph $\G$ is {\it edge-transitive} if $\Aut(\G)$ is transitive on the edge set $E(\G)$.

Let $\G$ be a graph and $N\leq \Aut(\G)$. The {\em quotient graph} $\G_N$ of $\G$ relative to the orbits
of $N$ is defined as the graph with vertices the orbits of $N$ on $V(\G)$ and with two orbits
adjacent if there is an edge in $\G$ between those two orbits.
In particular, for a normal subgroup $N$ of $\Aut(\G)$,
if $\G$ and $\G_N$ have the same valency, then $\G_N$ is a {\em normal quotient}
of $\G$, and if $\G$ has no proper normal quotient, then $\G$ is {\em basic}.
To study a symmetric graph $\G$, there is an extensive used strategy
consisting of two steps:
the first one is to investigate normal quotient graph
$\G_N$ for some normal subgroup $N$ of $\Aut(\G)$ and the second one is to reconstruct
the original graph $\G$ from the normal quotient $\G_N$ by using covering
techniques. This strategy was first laid out by Praeger (see \cite{P3}), and it is usually done by taking the normal subgroup $N$
as large as possible and then the graph $\G$ is reduced a basic graph.
In the literature, there are many works about basic graphs (see~\cite{FP,FZL} for example), while the works about the second step, that is, covers of graphs, are fewer.

An epimorphism $\pi: \widetilde{\G}\mapsto \G$ of graphs is called a {\em regular
covering projection} if $\Aut(\widetilde{\G})$ has a semiregular subgroup $K$ whose orbits in $V(\widetilde{\G})$ coincide with the {\em vertex fibres} $\pi^{-1}(v)$, $v\in V(\G)$, and whose arc and edge orbits coincide with the
{\em arc fibres} $\pi^{-1}((u,v))$ and the {\em edge fibres} $\pi^{-1}(\{u,v\})$, $\{u,v\}\in E(\G)$, respectively. In particular, we call the graph $\widetilde{\G}$ a {\em regular cover}
or a {\em $K$-cover} of the graph $\G$, and the group $K$ the {\em covering transformation group}.
If $K$ is dihedral, cyclic or elementary abelian, then $\widetilde{\G}$
is called a {\em dihedral}, {\em cyclic} or {\em elementary abelian cover} of $\G$, respectively.
An automorphism of $\widetilde{\G}$ is
said to be {\em fibre-preserving} if it maps a vertex fibre to a vertex fibre, and all
such fibre-preserving automorphisms form a group called the {\em fibre-preserving
group}, denoted by $F$. It is easy to see that $F=N_{\rm Aut(\widetilde{\G})}(K)$.
If $\widetilde{\G}$ is $F$-arc-transitive,
we say that $\widetilde{\G}$ is an {\em arc-transitive cover} or a
{\em symmetric cover} of $\G$. For an extensive treatment of regular cover,
one can see~\cite{MMP,27}.

Covering techniques have long been known as a powerful tool in algebraic and topological graph
theory. Application of these techniques has resulted in many
constructions and classifications of certain families of graphs with particular symmetry properties.
For example, by using covering techniques, Djokovi\' c~\cite{D}
constructed the first infinite family
of 5-arc transitive cubic graphs as covers of Tutte's 8-cage, and Biggs~\cite{B2} constructed some 5-arc-transitive cubic graphs as covers of
cubic graphs that are 4-arc-transitive but not 5-arc-transitive. Gross and Tucker~\cite{Gross} proved that every regular cover of a base graph
can be reconstructed as a voltage graph on the base graph.
Later, Malni\v c et al.~\cite{MMP} and Du et al.~\cite{Du2} developed these ideas
further in a systematic study of regular covering projections of a graph along which a group of automorphisms lifts.

Based on the approaches studied in~\cite{Du2,MMP},
many arc-transitive covers of symmetric graphs of
small orders and small valencies have been classified. For example,
Pan et al.~\cite{PHXD} studied arc-transitive cyclic covers of
some complete graphs of small orders. One may see~\cite{AHK,27} for other works.
Moreover, a new approach was proposed by Conder and Ma~\cite{Ma-covers,CM} by considering a presentation (quotient group) of a universal group, which can be obtained from Reidemeister-Schreier theory, and representation theory and other methods are applied when determining suitable quotients. As an application, arc-transitive abelian covers of the complete graph $K_4$, the complete bipartite graph $K_{3,3}$,
the 3-dimensional hypercube $Q_3$, the Petersen graph and the Heawood graph, were classified. Later, arc-transitive dihedral covers of these graphs were determined by Ma~\cite{Ma}.

For arc-transitive covers of infinite families of graphs,
Du et al. studied $2$-arc-transitive elementary abelian and cyclic covers of complete graphs $K_n$ in~\cite{Du1,Du4} and cyclic covers of $K_{n,n}-nK_2$ in~\cite{XD}. Recently, Pan et al.~\cite{PHL} determined arc-transitive cyclic covers of the complete bipartite graph $K_{p,p}$ of order $2p$ for a prime $p$. Compared with symmetric covers of graphs of small orders and valencies, there are only a few contributions on symmetric covers of infinite families of graphs.

Arc-transitive covers of non-simple graphs were also considered in literature.
For example, regular covers of the dipole $\Dip_k$ (a graph with two vertices and $k$ parallel edges) were extensively studied in~\cite{AHK,FZL,MMP,27,YF}. Such covers are called {\em Haar graphs},
and in particular, cyclic regular covers of dipoles are called {\em cyclic Haar graphs},
which can be regarded as a generalization of bipartite circulants and were studied in~\cite{HMP}
(also see~\cite{FL}). Construction of Haar graphs have aroused wide concern. Maru\v si\v c et al.~\cite{MMP} studied elementary abelian covers of the dipole $\Dip_p$
for a prime $p$. In particular, symmetric elementary abelian covers
and $\mz_p^2\times\mz_p$-covers for a prime $p$ of the dipole $\Dip_5$
were classified completely in~\cite{FZL} and~\cite{YF}, respectively.

Let $p$ be a prime. Pentavalent symmetric graphs of order $2p$ were classified by Cheng and Oxley
in~\cite{Cheng}, which are the complete graph $K_6$ of order 6 and a family of Cayley graphs $\mathcal{CD}_p$ with $p=5$ or $5\di (p-1)$ on dihedral groups (see Proposition~\ref{prop=2p}).
It has been shown that many pentavalent symmetric graphs are
regular covers of them, see~\cite{FZL,YF}. In this paper, we consider arc-transitive cyclic and dihedral covers of these graphs.
For $K_6$, the cyclic covers have been classified in~\cite{PHXD},
which should be the complete bipartite graph $K_{6,6}$ and the Icosahedron graph $\mathbf{I}_{12}$
(note that $\mathbf{I}_{12}$ is missed in \cite{PHXD}).
For $\mathcal{CD}_p$,
the cyclic covers consist of six infinite families of graphs, which are
Cayley graphs on generalized dihedral groups. In particular,
one family of graphs are cyclic Haar graphs and the other five families
are non-cyclic Haar graphs. What is more, the full automorphism groups
of them are determined. Arc-transitive dihedral covers
of $K_6$ and $\mathcal{CD}_p$ are also classified,
and there are only four sporadic graphs of order 24, 48, 60 and 120, respectively.
A similar work about cubic graphs was done by Zhou and Feng~\cite{ZF}.

Different from regular covers of graphs mentioned above,
the method to classify arc-transitive cyclic covers used in this paper is
related to the so called bi-Cayley graph.
A graph $\G$ is a {\em bi-Cayley graph} over some group $H$
if $\Aut(\G)$ has a semiregular subgroup isomorphic to $H$ having exactly two orbits on $V(\G)$. Clearly, a Haar graph is a bipartite bi-Cayley graph. Recently, Zhou and Feng~\cite{ZF1}
gave a depiction of the automorphisms of bi-Cayley graphs (see Section~4),
and based on this work, we classify the cyclic covers. In particular, all these covers are bi-Cayley graphs over some abelian groups.
Note that vertex-transitive bi-Cayley graphs of valency $3$ over abelian groups were determined in~\cite{ZF2}, while
the case for valency 5 is still elusive. Indeed, even for arc-transitive pentavalent bi-Cayley graphs over abelian groups,
it seems to be very difficult to give a classification, and one may see~\cite{AHK,FZL, YF} for partial works.

The paper is organized as follows. After this introductory section, in Section 2 we give some notation and preliminary results. In Section 3, several infinite families of connected pentavalent symmetric graphs are constructed as Cayley graphs
on generalized dihedral groups $\Dih(\mz_{mp^e}\times\mz_p)$, where $e,m$
are two positive integers and $p$ is a prime such that $(m,p)=1$.
In Section 4, it is proved that these Cayley graphs include all arc-transitive normal bipartite bi-Cayley graphs over $\mz_{mp^e}\times\mz_p$, and using this result, arc-transitive cyclic and dihedral covers of connected pentavalent symmetric graphs of order $2p$ are classified in Sections 5 and 6, respectively. In Section 7, the full automorphism groups of these covers are determined.

\section{Preliminaries}


In this section, we describe some preliminary results which will be used later.
The following result is important to investigate symmetric pentavalent graphs.

\begin{prop}\label{prop=stabilizer} {\rm \cite[Theorem 1.1]{GF}}
Let $\Gamma$ be a connected pentavalent $(G,s)$-transitive graph for some $G\leq {\rm Aut}(\Gamma)$ and $s\geq1$, and let $v\in V(\Gamma)$. Then one of the following holds:
\begin{enumerate}

\itemsep -1pt
\item [\rm (1)] $s=1$ and $G_v\cong \mathbb{Z}_5$, $D_{5}$ or $D_{10}$;
\item [\rm (2)] $s=2$ and $G_v\cong F_{20}$, $F_{20}\times \mathbb{Z}_2$, $A_5$ or $S_5$, where $F_{20}$ is the Frobenius group of order
$20$;
\item [\rm (3)] $s=3$ and $G_v\cong F_{20}\times \mathbb{Z}_4$, $A_4\times A_5$, $S_4\times S_5$ or $(A_4\times A_5)\rtimes \mathbb{Z}_2$ with $A_4\rtimes \mathbb{Z}_2=S_4$ and $A_5\rtimes \mathbb{Z}_2=S_5$;
\item [\rm (4)] $s=4$ and $G_v\cong {\rm ASL}(2,4)$, ${\rm AGL}(2,4)$, ${\rm A\Sigma L}(2,4)$ or ${\rm A\Gamma L}(2,4)$;
\item [\rm (5)] $s=5$ and $G_v\cong \mathbb{Z}^6_2\rtimes  {\rm \Gamma L}(2,4)$.
\end{enumerate}
\end{prop}

From~\cite[Theorem 9]{Lorimer}, we have the following proposition.

\begin{prop}\label{prop=atlesst3orbits}
Let $\Gamma$ be a connected $G$-arc-transitive graph of prime valency, and let $N$ be a normal subgroup of $G$.
If $N$ has at least three orbits, then it is semiregular on $V(\G)$ and the kernel of $G$ on the quotient graph $\G_N$. Furthermore, $\G_N$ is $G/N$-arc-transitive, and $\G$ is a regular cover of $\G_N$ with $N$ as the covering transformation group.
\end{prop}

Let $G$ and $E$ be two groups. We call an extension $E$ of $G$ by $N$ a
{\em central extension} of $G$ if $E$ has a central subgroup $N$ such that $E/N\cong G$,
and if further $E$ is perfect, that is, the derived group $E'=E$, we
call $E$ a {\em covering group} of $G$. Schur proved that for every non-abelian simple group $G$
there is a unique maximal covering group $M$ such that every covering group of $G$ is a factor
group of $M$ (see \cite[V \S23]{Huppert}). This group $M$ is called the {\em full covering group} of $G$, and the center of $M$
is the {\em Schur multiplier} of $G$, denoted by $\Mult(G)$.

\begin{lem}\label{lem=Mult}
Let $G$ be a group, and let $N$ be an abelian normal subgroup of $G$ such that
$G/N$ is a non-abelian simple group.
If $N$ is a proper subgroup of $C_G(N)$, then
$G=G'N$ and $G'\cap N\lesssim \Mult(G/N)$.
\end{lem}

\f {\bf Proof:} Since $N$ is a proper subgroup of $C_G(N)$, we have $1\neq C_G(N)/N \unlhd G/N$, forcing $C_G(N)/N=G/N$ because $G/N$ is simple. Thus $G=C_G(N)$ is a central extension of $G/N$ by $N$.
Since $G/N=(G/N)'=G'N/N\cong G'/(G'\cap N)$, we have $G=G'N$, and since $G'=(G'N)'=(G')'$, $G'$ is a
covering group of $G/N$. Hence $G'\cap N\lesssim \Mult(G/N)$.
\hfill$\blacksquare$

\medskip

Denote by $\soc(G)$ the {\em socle} of $G$, that is,
the product of all minimal normal subgroups of $G$.
A list of all proper primitive permutation groups of degree less than 1000 was given
by Dixon and Mortimer~\cite[Appendix B]{DM}, and based on the list, we have:

\begin{lem}\label{lem=primitive permutation group}
Let $G$ be a primitive permutation group on a set $\Omega$ and let $\a\in \Omega$,
where $|\Omega|\in \{2,4,6,8,12,16,24,72,144,288,576\}$.
If $G_{\a}$ is solvable, then either $G\lesssim {\rm AGL}(n,2)$ and
$|\Omega|=2^n$ with $1\leq n\leq 4$, or $\soc(G)\cong \PSL(2,p)$,
$\PSL(3,3)$ or $\PSL(2,q)\times\PSL(2,q)$ with $|\Omega|=p+1$, $144$ or $(q+1)^2$ respectively, where $p\in \{5,7,11,23,71\}$ and $q\in \{11,23\}$.
\end{lem}

\f {\bf Proof:} If $|\Omega|=2$ or 4, then $G\leq S_{2}\cong {\rm AGL}(1,2)$
or $G\leq S_4\cong {\rm AGL}(2,2)$, respectively. Let $|\Omega|\geq 6$ and write $N:=\soc(G)$. Then $N\unlhd G$ and $N_{\a}\unlhd G_{\a}$.
Since $G_{\a}$ is solvable, $N_{\a}$ is solvable.
By~\cite[Appendix B, Tables~B.2 and B.3]{DM}, $G$ is an affine group, $N\cong A_{|\Omega|}$, or
$G$ is isomorphic to one group listed in~\cite[Tables~B.2 and B.3]{DM}.
If $G$ is affine, then $|\Omega|$ is a prime power and thus $|\Omega|=2^n$ with $n=3$ or 4.
By \cite[Theorem~4.1A (a)]{DM}, we have $G\lesssim {\rm AGL}(n,2)$.
If $N\cong A_{|\Omega|}$ then $N_{\a}\cong A_{|\Omega|-1}$, which is insolvable because
$|\Omega|-1\geq 5$, a contradiction.
In what follows we assume that $G$ is isomorphic to one group listed in~\cite[Tables B.2 and B.3]{DM}.
Note that all groups in the tables are collected into cohorts
and all groups in a cohort have the same socle.

Assume that $|\Omega|=144$. By~\cite[Table B.4, pp.324]{DM}, there
are one cohort of type $C$, two cohorts of type $H$ and four cohorts of type $I$
(see~\cite[Table~B.1, pp.306]{DM} for types of cohorts of primitive groups) of primitive groups of degree 144.
For the cohort of type $C$, by~\cite[Table~B.2, pp.314]{DM}, $N\cong \PSL(3,3)$
and $N_{\a}\cong \mz_{13}\rtimes\mz_3$.
For the two cohorts of type $H$, by~\cite[Table~B.2, pp.321]{DM},
they have the same socle $N\cong {\rm M}_{12}$ and $N_{\a}\cong \PSL(2,11)$.
For the four cohorts of type $I$, by~\cite[Table~B.3, pp.323]{DM},
$N\cong A_{12}\times A_{12}$, $\PSL(2,11)\times\PSL(2,11)$,
${\rm M}_{11}\times {\rm M}_{11}$ or ${\rm M}_{12}\times {\rm M}_{12}$
and $N_{\a}\cong A_{11}\times A_{11}$, $(\mz_{11}\rtimes\mz_{5})\times(\mz_{11}\rtimes\mz_{5})$,
${\rm M}_{10}\times {\rm M}_{10}$ or ${\rm M}_{11}\times {\rm M}_{11}$, respectively.
Since $N_{\a}$ is solvable, we have $N\cong \PSL(3,3)$ or $\PSL(2,11)\times\PSL(2,11)$.

For $|\Omega|\in \{6,8,12,16,24,72,288,576\}$,
by~\cite[Tables~B.2, B.3 and B.4]{DM}, a similar argument to the above paragraph implies that either $N\cong  \PSL(2,23)\times \PSL(2,23)$ with degree $23^2=576$ and $N_{\a}\cong (\mz_{23}\rtimes\mz_{11})\times(\mz_{23}\rtimes\mz_{11})$, or $N\cong\PSL(2,p)$ with degree
$p+1$ and $N_{\a}\cong \mz_p\rtimes\mz_{\frac{p-1}{2}}$ where $p\in \{5,7, 11, 23, 71\}$.
\hfill$\blacksquare$

\section{Graph constructions as Cayley graphs}

Let $G$ be a finite group and $S$ a subset of $G$ with $1\not\in S$
and $S^{-1}=S$. The {\em Cayley graph} $\G=\Cay(G,S)$ on $G$ with
respect to $S$ is defined to have vertex set $V(\G)=G$ and edge set
$E(\G)=\{\{g,sg\}\ |\ g\in G,s\in S\}$. It is well-known that
$\Aut(\G)$ contains the right regular representation $R(G)$ of $G$,
the acting group of $G$ by right multiplication, and $\G$ is
connected if and only if $G=\langle S \rangle$, that is, $S$
generates $G$. By Godsil~\cite{G},
$N_{{\rm Aut}(\G)}(R(G))=R(G)\rtimes\Aut(G,S)$, where $\Aut(G,S)=\{\a\in \Aut(G) \ |\ S^\a=S\}$. A Cayley graph
$\G=\Cay(G,S)$ is said to be {\em normal} if $R(G)$ is normal in $\Aut(\G)$, and in
this case, $\Aut(\G)=R(G)\rtimes\Aut(G,S)$.

For an abelian group $H$,
the {\em generalized dihedral group} ${\rm Dih}(H)$ is the semidirect product $H\rtimes \mz_2$,
where the unique involution in $\mz_2$ maps each element of $H$ to its inverse. In
particular, if $H$ is cyclic, then ${\rm Dih}(H)$ is a dihedral
group. In this section, we introduce several infinite families of connected
pentavalent symmetric graphs which are constructed as Cayley graphs
on generalized dihedral groups.

\begin{exam}{\rm \label{exam=250}
Let $\Dih(\mz_5^3)=\lg a,b,c,h~|~a^5=b^5=c^5=h^2=[a,b]=[a,c]=[b,c]=1,a^h=a^{-1},
b^h=b^{-1},c^h=c^{-1}\rg$, and define
$$\mathcal{CGD}_{5^3}=\Cay(\Dih(\mz_5^3), \{h,ah,bh,ch,a^{-1}b^{-1}c^{-1}h\}).$$
By~\cite[Theorem~1.1]{YF}, $\Aut(\mathcal{CGD}_{5^3})\cong \Dih(\mz_5^3)\rtimes S_5$
and $\mathcal{CGD}_{5^3}$ is the unique connected pentavalent graph of order $250$ up to isomorphism.

}
\end{exam}

Let $m$ be a positive integer. Consider the following equation in $\mz_m$

\parbox{8cm}{
\begin{eqnarray*}
&&x^4+x^3+x^2+x+1=0.
\end{eqnarray*}}\hfill
\parbox{1cm}{
\begin{eqnarray}
\label{eq-r}
\end{eqnarray}}

\f In view of~\cite[Lemma~3.3]{FL}, we have the following proposition.

\begin{prop}\label{prop=x5=1}
Eq~{\rm (\ref{eq-r})} has a solution $r$ in $\mz_m$ if and only if
$(r,m)\in \{(0,1),(1,5)\}$ or
$m=5^tp_1^{e_1}p_2^{e_2}\cdots p_s^{e_s}$
and $r$ is an element in $\mz_m^*$ of order $5$, where
$t\leq 1$, $s\geq1$, $e_i\geq1$ and $p_i$'s are distinct
primes such that $5\ |\ (p_i-1)$.
\end{prop}

The following infinite
family of Cayley graphs was first constructed in~\cite{KKO}.

\begin{exam}{\rm \label{exam=n}
Let $m>1$ be an integer such that Eq~(\ref{eq-r}) has a
solution $r$ in $\mz_m$. Then $m=5$, 11 or $m\geq31$. Let
$$\mathcal{CD}_{m}=\Cay(D_{m},\{b,ab,a^{r+1}b,a^{r^2+r+1}b,a^{r^3+r^2+r+1}b\})$$
be a Cayley graph on the dihedral group
$D_{m}=\lg a,b~|~a^{n}=b^2=1, a^b=a^{-1}\rg$.
For $m=5$ or 11, by~\cite{Cheng},
$\Aut(\mathcal{CD}_m)\cong (S_5\times S_5)\rtimes\mz_2$
or $\PGL(2,11)$, respectively. In particular, $\mathcal{CD}_{5}\cong K_{5,5}$.
For $m\geq31$, by~\cite[Theorem~B and Proposition~4.1]{KKO},
$\Aut(\mathcal{CD}_{m})\cong D_{m}\rtimes\mz_5$, and obviously,
if $m$ has a prime divisor $p$ with $p<m$, then
$\Aut(\mathcal{CD}_{m})$ has a normal subgroup $\mz_{m/p}$, and by Proposition~\ref{prop=atlesst3orbits}, $\mathcal{CD}_{m}$ is a symmetric $\mz_{m/p}$-cover of a connected pentavalent symmetric graph of order $2p$.

}\end{exam}

By~\cite{Cheng}, we have the following proposition.

\begin{prop}\label{prop=2p}
Let $\G$ be a connected  pentavalent symmetric graph of order $2p$ for a prime $p$.
Then $\G\cong K_6$ or $\mathcal{CD}_{p}$ with $p=5$
or $5\ |\ (p-1)$.
\end{prop}

In the remaining part of this section, we construct five infinite families of Cayley graphs on some generalized dihedral groups, and for convenience, we always assume that $G=\Dih(\mz_m\times\mz_{p^e}\times\mz_p)=\lg a,b,c,h~|~a^m=b^{p^e}=c^p=h^2=[a,b]=[a,c]=[b,c]=1,
a^h=a^{-1},b^{h}=b^{-1},c^{h}=c^{-1}\rg$ and $r$ is a solution of Eq~(\ref{eq-r}) in $\mz_m$, that is, $r^4+r^3+r^2+r+1=0~(\mod m)$.
By Proposition~\ref{prop=x5=1}, $m$ is odd and $5^2\nmid m$.

\begin{exam}{\rm \label{exam=n*p}
Assume that $e\geq2$ and $p$ is a prime such that $(m,p)=1$
and $5\ |\ (p-1)$. Let $\lambda$ be an element of order $5$ in $\mz_{p^e}^*$. Then $\lambda$ is a solution of Eq~(\ref{eq-r})
in $\mz_{p^e}$. Set
$$\begin{array}{lll}
&T_1(r,\lambda)=\{h, hab, ha^{r+1}b^{\lambda+1}c, ha^{r^2+r+1}b^{\lambda^2+\lambda+1}c^{\lambda^4+\lambda+1},
ha^{r^3+r^2+r+1}b^{\lambda^3+\lambda^2+\lambda+1}c\},\\
&T_2(r,\lambda)=\{h, hab, ha^{r+1}b^{\lambda+1}c, ha^{r^2+r+1}b^{\lambda^2+\lambda+1}c^{\lambda^3+\lambda+1},
ha^{r^3+r^2+r+1}b^{\lambda^3+\lambda^2+\lambda+1}c^{\lambda}\},\\
&T_3(r,\lambda)=\{h,hab, ha^{r+1}b^{\lambda+1}c, ha^{r^2+r+1}b^{\lambda^2+\lambda+1}c^{\lambda^2+\lambda+1},
ha^{r^3+r^2+r+1}b^{\lambda^3+\lambda^2+\lambda+1}c^{\lambda^2}\}.
\end{array}$$
It is easy to see that each of $T_1(r,\lambda)$, $T_2(r,\lambda)$ and $T_3(r,\lambda)$ generates $G$. Define
$$\mathcal{CGD}_{mp^e\times p}^{i}=\Cay(G, T_i(r,\lambda)),~i=1,2,3.$$
The maps $a\mapsto a^r$, $b\mapsto b^{\lambda}c$,
$c\mapsto c^{\lambda^4}$, $h\mapsto hab$;
$a\mapsto a^r$, $b\mapsto b^{\lambda}c$,
$c\mapsto c^{\lambda^3}$, $h\mapsto hab$;
$a\mapsto a^r$, $b\mapsto b^{\lambda}c$,
$c\mapsto c^{\lambda^2}$, $h\mapsto hab$
induce three automorphisms of order $5$ of $G$ ,
denoted by $\a_1$, $\a_2$ and $\a_3$ respectively,
and $\a_i$ fixes the set $T_i(r,\lambda)$ and permutes
its five elements cyclicly. It follows that for each $i=1,2,3$, $\a_i\in\Aut(G,T_i(r,\lambda))$ and
$\lg R(G),\a_i\rg\cong G\rtimes \mz_5$, which is arc-transitive on $\mathcal{CGD}_{mp^e\times p}^{i}$.
}
\end{exam}

\begin{lem}\label{lem=exam4.1}
The graphs $\mathcal{CGD}_{mp^e\times p}^{1}$, $\mathcal{CGD}_{mp^e\times p}^{2}$ and $\mathcal{CGD}_{mp^e\times p}^{3}$ are
not isomorphic to each other.
\end{lem}

\f {\bf Proof:} Let $\G=\mathcal{CGD}_{mp^e\times p}^{i}$
and $A=\Aut(\G)$, where $i=1,2$ or $3$.
Then $A=R(G)\rtimes \lg \a_i\rg
\cong(\mz_{mp^e}\times\mz_p)\rtimes\mz_2)\rtimes \mz_5$ because $(m,p)=1$.
First, we claim that all regular subgroups of $A$ isomorphic to $G$
are conjugate. Recall that $p>5$ and $5^2\nmid m$.
If $5\nmid m$ then $G$ is a Hall $5'$-subgroup of $A$
and thus all the regular subgroups isomorphic to $G$
are conjugate in $A$. If $5\di m$ then $A$ has a
characteristic Hall $\{2,5\}'$-subgroup, say $P$. It follows $|P|=mp^{e+1}/5$ and $|A/P|=50$.
Clearly, $P$ has more than two orbits on $V(\G)$
and $P$ is contained in every regular subgroup
of $A$. By Proposition~\ref{prop=atlesst3orbits},
$\G_P$ is a connected pentavalent $A/P$-arc-transitive graph of order 10
and by Proposition~\ref{prop=2p}, $\G_P\cong K_{5,5}$.
Let $G_1$ and $G_2$
be any two regular subgroups such that $G_1\cong G_2\cong G$.
As $G$ is a generalized dihedral group, we have
$G_1/P\cong G_2/P\cong D_5$.
Note that $A/P$ is an arc-transitive subgroup of
$\Aut(K_{5,5})$ of order $50$, and by MAGMA~\cite{BCP},
all subgroups of $A/P$ isomorphic to $D_5$ are conjugate in $A/P$.
It implies that there is an element $gP\in A/P$ such that
$(G_1/P)^{gP}=G_1^g/P=G_2/P$, where $g\in A$. Hence $G_1^g=G_2$, as claimed.

Now we are ready to prove that the graphs $\mathcal{CGD}_{mp^e\times p}^{1}$,
$\mathcal{CGD}_{mp^e\times p}^{2}$ and $\mathcal{CGD}_{mp^e\times p}^{3}$ are
not isomorphic to each other.
Suppose that $\mathcal{CGD}_{mp^e\times p}^{1}\cong \mathcal{CGD}_{mp^e\times p}^{2}$.
Since all regular subgroups isomorphic to $G$ in $\Aut(\mathcal{CGD}_{mp^e\times p}^{1})$
are conjugate, by Babai~\cite[Lemma~3.1]{Babai}, there exists
an automorphism  $\b\in \Aut(G)$ such that $T_1(r,\lambda)^{\b}=T_2(r,\lambda)$, that is,
$$\begin{array}{ll}
&~~~\{h, hab, ha^{r+1}b^{\lambda+1}c, ha^{r^2+r+1}b^{\lambda^2+\lambda+1}c^{\lambda^4+\lambda+1},
ha^{r^3+r^2+r+1}b^{\lambda^3+\lambda^2+\lambda+1}c\}^\b\\
&=\{h, hab, ha^{r+1}b^{\lambda+1}c, ha^{r^2+r+1}b^{\lambda^2+\lambda+1}c^{\lambda^3+\lambda+1},
ha^{r^3+r^2+r+1}b^{\lambda^3+\lambda^2+\lambda+1}c^{\lambda}\}.
\end{array}$$

By Example~\ref{exam=n*p}, $\a_2\in\Aut(G)$ permutes the elements in $T_2(r,\lambda)$ cyclicly, and thus we may assume that $h^{\b}=h$, which implies the following equation:

\parbox{8cm}{
\begin{eqnarray*}
&&~~~\{ab, a^{r+1}b^{\lambda+1}c, a^{r^2+r+1}b^{\lambda^2+\lambda+1}c^{\lambda^4+\lambda+1},
a^{r^3+r^2+r+1}b^{\lambda^3+\lambda^2+\lambda+1}c\}^\b\\
&&=\{ab, a^{r+1}b^{\lambda+1}c, a^{r^2+r+1}b^{\lambda^2+\lambda+1}c^{\lambda^3+\lambda+1},
a^{r^3+r^2+r+1}b^{\lambda^3+\lambda^2+\lambda+1}c^{\lambda}\}.
\end{eqnarray*}}\hfill
\parbox{1cm}{
\begin{eqnarray}
\label{eq4.2}
\end{eqnarray}}

\f By assumption in Example~\ref{exam=n*p}, $(m,p)=1$, and since $m$ is odd, $(m,2p)=1$. Thus both $\mz_m=\lg a\rg$ and $\mz_{p^e}\times\mz_p=\lg b,c\rg$ are characteristic in $G$, and from Eq~(\ref{eq4.2}) we have

\parbox{1cm}{
\begin{eqnarray*}
&\{b, b^{\lambda+1}c, b^{\lambda^2+\lambda+1}c^{\lambda^4+\lambda+1},
b^{\lambda^3+\lambda^2+\lambda+1}c\}^\b=\{b, b^{\lambda+1}c, b^{\lambda^2+\lambda+1}c^{\lambda^3+\lambda+1},
b^{\lambda^3+\lambda^2+\lambda+1}c^{\lambda}\}.
\end{eqnarray*}}\hfill
\parbox{1cm}{
\begin{eqnarray}
\label{eq4.2*}
\end{eqnarray}}

\f It follows that $b^\b=b^sc^t$, where $(s,t)=(1,0),(\lambda+1,1), (\lambda^2+\lambda+1,\lambda^3+\lambda+1)$ or $(\lambda^3+\lambda^2+\lambda+1,\lambda)$. Furthermore, we have
$$\begin{array}{l}
(b\cdot b^{\lambda+1}c\cdot b^{\lambda^2+\lambda+1}c^{\lambda^4+\lambda+1}\cdot
b^{\lambda^3+\lambda^2+\lambda+1}c)^{\b}=b\cdot b^{\lambda+1}c\cdot b^{\lambda^2+\lambda+1}c^{\lambda^3+\lambda+1}\cdot
b^{\lambda^3+\lambda^2+\lambda+1}c^{\lambda},
\end{array}$$
that is, $(b^\b b^{-1})^{-\lambda^4+\lambda^2+2\lambda+3}
=(c^{\b})^{-\lambda^4-\lambda-3}c^{\lambda^3+2\lambda+2}$. In particular,
$(b^\b b^{-1})^{(-\lambda^4+\lambda^2+2\lambda+3)p}=1$.

Suppose that $-\lambda^4+\lambda^2+2\lambda+3=0$ in $\mz_p$. Since
$\lambda^4+\lambda^3+\lambda^2+\lambda+1=0$ in $\mz_{p^e}$, it is also true in $\mz_p$ and hence $\lambda^5=1$. Thus $\lambda^4=\lambda^2+2\lambda+3$ and $\lambda^3=\lambda\cdot \lambda^2=\lambda(\lambda^4-2\lambda-3)=-2\lambda^2-3\lambda+1$, which implies that  $0=\lambda^4+\lambda^3+\lambda^2+\lambda+1=5$ in $\mz_p$, contrary to the assumption that $5\ |\ (p-1)$. It follows that $-\lambda^4+\lambda^2+2\lambda+3\not=0$ in $\mz_p$ and  $(b^\b b^{-1})^p=1$.

Suppose that $(s,t)\not=(1,0)$. Then $b^\b b^{-1}=b^{s-1}c^t$ with $s-1=\lambda,\lambda^2+\lambda$ or $\lambda^3+\lambda^2+\lambda$. Since $\lambda^4+\lambda^3+\lambda^2+\lambda+1=0$ in $\mz_p$, one may easily show that $(s-1,p)=1$. This implies that $b^\b b^{-1}=b^{s-1}c^t$ has order $p^e$, and since $e\geq 2$, we have  $(b^\b b^{-1})^p\not=1$, a contradiction. It follows that $(s,t)=(1,0)$, that is, $b^\b=b$, and by Eq~(\ref{eq4.2*}), $\{c, c^{\lambda^4+\lambda+1}, c\}^\b=\{c, c^{\lambda^3+\lambda+1}, c^{\lambda}\}$,
which is impossible because any two elements in $\{c, c^{\lambda^3+\lambda+1}, c^{\lambda}\}$
are not equal.
Hence $\mathcal{CGD}_{mp^e\times p}^{1}\ncong \mathcal{CGD}_{mp^e\times p}^{2}$.
Similarly, on may check that $\mathcal{CGD}_{mp^e\times p}^{1}\ncong \mathcal{CGD}_{mp^e\times p}^{3}$ and $\mathcal{CGD}_{mp^e\times p}^{2}\ncong \mathcal{CGD}_{mp^e\times p}^{3}$.\hfill$\blacksquare$

\begin{exam}{\rm \label{exam=n*p*p-1} Let $p$ be a prime such that $p=5$ or $5\ |\ (p\pm1)$.
Assume that $e=1$ and $(m,p)=1$. Then $G=\Dih(\mz_m\times\mz_p\times\mz_p)$.
For $p=5$, let $\lambda=0$, and for $5\ |\ (p\pm1)$,
let $\lambda\in\mz_p$ satisfying the equation $\lambda^2=5$ in $\mz_p$. Set
$$S(r,\lambda)=\{h,hab,ha^{r+1}c,ha^{r^2+r+1}b^{-2^{-1}(1+\lambda)}c^{2^{-1}(1+\lambda)},
ha^{r^3+r^2+r+1}b^{-2^{-1}(1+\lambda)}c\}.$$
It is easy to see that $S(r,\lambda)$ generates $G$.
Define
$$\mathcal{CGD}_{mp\times p}^{4}=\Cay(G,S(r,\lambda)).$$
The map $a\mapsto a^r$, $b\mapsto b^{-1}c$, $c\mapsto b^{-2^{-1}(3+\lambda)}c^{2^{-1}(1+\lambda)}$ and $h\mapsto hab$
induces an automorphism of the group $G$, denoted by $\a_4$, which permutes
the elements in $S(r,\lambda)$ cyclicly. Then $\a_4\in \Aut(G,S(r,\lambda))$
and $\lg R(G),\a_4\rg\cong G\rtimes \mz_5$
acts arc-transitive on $\mathcal{CGD}_{mp\times p}^{4}$.
Moreover, for $m=1$ or $5$, we have $r=0$ or $1$ respectively,
and the map $a\mapsto a^{-1}$, $b\mapsto b^{-2^{-1}(1+\lambda)}c$,
$c\mapsto b^{-2^{-1}(1+\lambda)}c^{2^{-1}(1+\lambda)}$, $h\mapsto h$ induces an automorphism $\b$ of $G$. It is easy to check that $\b\in \Aut(G,S(r,\lambda))$
and $\lg \a_4,\b\rg\cong D_5$.
In particular, by~\cite[Theorem~6.1]{FZL}, $\mathcal{CGD}_{5\times 5}^{4}$
is the unique connected pentavalent symmetric graph of order $50$ up to isomorphism. }
\end{exam}

\begin{exam}{\rm \label{exam=n*p*p-2}
Assume that $e=1$ and $p$ is a prime such that $(m,p)=1$ and $5\ |\ (p-1)$.
By~\cite[Case 2, page 14]{YF}, $x^{4}+10x^2+5=0$ has a root $\lambda$ in $\mz_p$.
Denote by $S(r,\lambda)$ the set
$$\{h,hab,ha^{r+1}c,ha^{r^2+r+1}b^{8^{-1}(\lambda^3-\lambda^2+7\lambda+1)}c^{2^{-1}(\lambda+1)},ha^{r^3+r^2+r+1}b^{-8^{-1}(\lambda^3+\lambda^2+7\lambda-1)}c^{8^{-1}(\lambda^3+\lambda^2+11\lambda+3)}\}.
$$
It is easy to check that $S(r,\lambda)$ generates $G$. Define
$$\mathcal{CGD}_{mp\times p}^{5}=\Cay(G,S(r,\lambda)).$$
The map $a\mapsto a^{r}$, $b\mapsto b^{-1}c$, $c\mapsto b^{8^{-1}(\lambda^3-\lambda^2+7\lambda-7)}
c^{2^{-1}(\lambda+1)}$ and $h\mapsto hab$
induces an automorphism of the group $G$, denoted by $\a_5$, which permutes the elements in
$S(r,\lambda)$ cyclicly. Then $\a_5\in \Aut(G,S(r,\lambda))$ and $\lg R(G),\a_5\rg\cong G\rtimes \mz_5$
acts arc-transitive on $\mathcal{CGD}_{mp\times p}^{5}$.
}\end{exam}

\begin{lem}\label{lem=examples}
Let $p$ be a prime such that $p=5$ or $5\di (p-1)$. Then for each $1\leq i\leq 5$, $\mathcal{CGD}_{mp^e\times p}^{i}$ is a connected symmetric cyclic cover of a connected
pentavalent symmetric graph of order $2p$.
\end{lem}
\f {\bf Proof:} Let $\G_i=\mathcal{CGD}_{mp^e\times p}^{i}$ for $i=1,2,3,4,5$. By definition of $\G_i$, $(m,p)=1$. Furthermore,
$|V(\G_i)|=2mp^{e+1}$ and $\Aut(\G_i)$ contains an arc-transitive subgroup $R(G)\rtimes\lg \a_i\rg$, where $\a_i$ for $1\leq i\leq 5$ are defined in Examples~\ref{exam=n*p}, \ref{exam=n*p*p-1} and \ref{exam=n*p*p-2}:
$$\begin{array}{ll}
&\a_i:~a\mapsto a^r,~b\mapsto b^{\lambda}c,~c\mapsto c^{\lambda^{5-i}},~h\mapsto hab,~i=1,2,3;\\
&\a_4:~a\mapsto a^r,~b\mapsto b^{-1}c,~c\mapsto b^{-2^{-1}(3+\lambda)}c^{2^{-1}(1+\lambda)},~h\mapsto hab;\\
&\a_5:~a\mapsto a^{r},~b\mapsto b^{-1}c,~c\mapsto b^{8^{-1}(\lambda^3-\lambda^2+7\lambda-7)}
c^{2^{-1}(\lambda+1)},~h\mapsto hab.
\end{array}$$

To prove that
$\G_i$ is a connected symmetric cyclic cover of a connected pentavalent symmetric graph of order $2p$, by Proposition~\ref{prop=atlesst3orbits}, it suffices to prove that $R(G)\rtimes\lg \a_i\rg$ has a normal cyclic subgroup $N_i\cong\mz_{mp^e}$ for each $1\leq i\leq 5$.

Assume $i=1$. Set $N_1=\lg R(a),R(b^{5}c^{3\lambda^4+2\lambda^2-\lambda+1})\rg$. Then $N_1\cong\mz_{mp^e}$. By Example~\ref{exam=n*p}, $\lambda^5=1$ and $\lambda^4+\lambda^3+\lambda^2+\lambda+1=0$ in $\mz_{p^e}$.  Since $a^{\a_1}=a^r$, $b^{\a_1}=b^{\lambda}c$ and $c^{\a_1}=c^{\lambda^4}$, we have
$$(b^{5}c^{3\lambda^4+2\lambda^2-\lambda+1})^{\a_1}
=b^{5\lambda}c^{\lambda^4+3\lambda^3+2\lambda+4}
=b^{5\lambda}c^{2\lambda^3-\lambda^2+\lambda+3}
=(b^{5}c^{3\lambda^4+2\lambda^2-\lambda+1})^{\lambda},$$
yielding that $\lg a,b^{5}c^{3\lambda^4+2\lambda^2-\lambda+1}\rg^{\a_1}=\lg a,b^{5}c^{3\lambda^4+2\lambda^2-\lambda+1}\rg$.
Since $R(x)^{\a_1}=R(x^{\a_1})$ for any $x\in G$,
we have $\lg R(a),R(b^{5}c^{3\lambda^4+2\lambda^2-\lambda+1})\rg^{\a_1}=N_1^{\a_1}=N_1$, and since $N_1\unlhd R(G)$, we have $N_1\unlhd \lg R(G), \a_1\rg$, as required.

Assume $i=2$ or $3$. Similar to the case $i=1$, one may check that $N_i\unlhd \lg R(G),\a_i\rg$ for $N_2=\lg R(a),R(b^{-5}c^{2\lambda^3+4\lambda^2+\lambda+3})\rg$ and $N_3=\lg R(a),R(b^{-5}c^{4\lambda^3+3\lambda^2+2\lambda+1})\rg$, as required.

Assume $i=4$. By assumption, $p=5$ or $5\di (p-1)$. Note that $\G_4\cong \mathcal{CGD}_{mp\times p}^{4}$. For $p=5$, by Example~\ref{exam=n*p*p-1}, $\lambda=0$ and thus $a^{\a_4}=a^r$, $b^{\a_4}=b^{-1}c$, $c^{\a_4}=bc^3$.
Set $N_4=\lg R(a),R(b^2c^4)\rg$. Then $(b^2c^4)^{\a_4}=b^2c^4$, $N_4\unlhd \lg R(G), \a_4\rg$ and $N_4\cong\mz_{5m}$, as required. For $5\di (p-1)$,
by Example~\ref{exam=n*p*p-2} the equation $x^{4}+10x^2+5=0$ has a root in $\mz_p$, say $t$,
that is, $t^{4}+10t^2+5=0$. Then $[2^{-1}(t^2+5)]^2=5=t^2(-t^2-10)$, and
thus $-t^2-10=5t^{-2}=[2^{-1}(t^2+5)]^2\cdot t^{-2}=[2^{-1}(t+5t^{-1})]^2$.
Since $\lambda^2=5$ (see Example~\ref{exam=n*p*p-1}), we have $\lambda=\pm 2^{-1}(t^2+5)$
and $2\lambda-5=t^2$ or $[2^{-1}(t+5t^{-1})]^2(=-t^2-10)$.
If $2\lambda-5=t^2$, set $N_4=\lg R(a),R(b^{t+1}c^{\lambda-3})\rg$. Since $a^{\a_4}=a^r$, $b^{\a_4}=b^{-1}c$ and
$c^{\a_4}=b^{-2^{-1}(3+\lambda)}c^{2^{-1}(1+\lambda)}$, we have
$$(b^{t+1}c^{\lambda-3})^{\a_4}=b^{1-t}c^{2+t-\lambda}=
(b^{t+1}c^{\lambda-3})^{-4^{-1}[t(\lambda+3)-\lambda+1]}\in \lg b^{t+1}c^{\lambda-3}\rg,$$
implying that $\lg a,b^{t+1}c^{\lambda-3}\rg^{\a_4}=\lg a,b^{t+1}c^{\lambda-3}\rg$
and $\lg R(a),R(b^{t+1}c^{\lambda-3})\rg^{\a_4}=N_4^{\a_4}=N_4$.
Since $N_4\unlhd R(G)$, we have $N_4\unlhd \lg R(G), \a_4\rg$, as required.
If $2\lambda-5=[2^{-1}(t+5t^{-1})]^2$, set $N_4=\lg R(a),R(b^{2^{-1}(t+5t^{-1})+1}c^{\lambda-3})\rg$ and
the above argument implies that
$$(b^{2^{-1}(t+5t^{-1})+1}c^{\lambda-3})^{\a_4}=(b^{2^{-1}(t+5t^{-1})+1}c^{\lambda-3})^{-4^{-1}[2^{-1}(t+5t^{-1})(\lambda+3)-\lambda+1]}\in \lg b^{2^{-1}(t+5t^{-1})+1}c^{\lambda-3}\rg.$$
Thus $\lg R(a),R(b^{2^{-1}(t+t^{-1}+1}c^{\lambda-3})\rg^{\a_4}=N_4^{\a_4}=N_4$
and $N_4\unlhd \lg R(G), \a_4\rg$, as required.

Assume $i=5$. By Example~\ref{exam=n*p*p-2}, we have
$5\di (p-1)$ and $\lambda^4+10\lambda^2+5=0$ in $\mz_p$.
Then $[2(\lambda^2+5)^{-1}]^2=5^{-1}$.
Set $N_5=\lg R(a),R(b^{t(\lambda^3+10\lambda+5)-(\lambda+3)}c^{4})\rg$,
where $t=2(\lambda^2+5)^{-1}$.
Then $N_5\cong\mz_{mp}$. Note that
$\lambda^4=-10\lambda^2-5$, $\lambda^5=-10\lambda^3-5\lambda$ and $\lambda^6=95\lambda^2+50$.
Since
$a^{\a_5}=a^{r}$, $b^{\a_5}=b^{-1}c$ and $c^{\a_5}=b^{8^{-1}(\lambda^3-\lambda^2+7\lambda-7)}
c^{2^{-1}(\lambda+1)}$, we have
$$\begin{array}{ll}
&(b^{t(\lambda^3+10\lambda+5)-(\lambda+3)}c^{4})^{\a_5}
=b^{-t(\lambda^3+10\lambda+5)+2^{-1}(\lambda^3-\lambda^2+9\lambda-1)}
c^{t(\lambda^3+10\lambda+5)+\lambda-1}\\
&=(b^{t(\lambda^3+10\lambda+5)-(\lambda+3)}c^{4})^{4^{-1}[t
(\lambda^3+10\lambda+5)+\lambda-1]}\in\lg b^{t(\lambda^3+10\lambda+5)-(\lambda+3)}c^{4}\rg.
\end{array}$$
Thus $N_5^{\a_5}=N_5$ and $N_5\unlhd \lg R(G), \a_5\rg$, as required.  \hfill$\blacksquare$

\section{Pentavalent symmetric bi-Cayley graphs over abelian groups}

Given a group $H$, let $R$, $L$ and $S$ be three subsets of $H$ such that
$R^{-1}=R$, $L^{-1}=L$, and $1\notin R\cup L$.
The {\em bi-Cayley graph} over $H$ relative to the triple $(R, L, S)$, denoted by
$\BiCay(H, R, L, S)$, is the graph
having vertex set $\{h_0~|~h\in H\}\cup \{h_1~|~h\in H\}$ and edge set $\{\{h_0, g_0\}~|~gh^{-1}\in R\}\cup \{\{h_1, g_1\}~|~gh^{-1}\in L\}\cup \{\{h_0, g_1\}~|~ gh^{-1}\in S\}$.
For a bi-Cayley graph $\G=\BiCay(H, R, L, S)$, it is easy to see that $R(H)$ can be regarded as a semiregular subgroup of $\Aut(\G)$, which acts on $V(\G)$ by the rule $h_i^{R(g)}=(hg)_i,~i=0,1,~h,g\in H$. If $R(H)$ is normal in $\Aut(\G)$, then
$\G$ is a {\em normal bi-Cayley graph} over $H$.

Let $\G=\BiCay(H, \emptyset, \emptyset, S)$ be a connected bi-Cayley graph over an abelian group $H$. Then $\G$ is bipartite. By~\cite[Lemma~3.1]{ZF1}, we may always assume that $1\in S$.
Moreover, $\G\cong\BiCay(H, \emptyset, \emptyset, S^{\a})$ for $\a\in\Aut(H)$, and $H=\lg S\rg$.
Since $H$ is abelian, there is an automorphism of $H$ of order 2, denoted by $\g$,
induced by $g\mapsto g^{-1},~\forall g\in H$.
For $\a\in \Aut(H)$ and $x\in H$, define

$$\begin{array}{ll}
&\d_{\g,1,1}:~h_0\mapsto (h^{-1})_1,~h_1\mapsto (h^{-1})_0,~\forall h\in H;\\
&\sigma_{\a,x}:~h_0\mapsto (h^{\a})_0,~h_1\mapsto (xh^{\a})_1,~\forall h\in H.
\end{array}$$
\f Set
$${\rm F}=\{\sigma_{\a,x}~|~\a\in\Aut(H),~S^{\a}=x^{-1}S\}.$$
\f Then $\d_{\g,1,1}\in\Aut(\G)$ and ${\rm F}\leq \Aut(\G)_{1_0}$
(see~\cite[Lemma~3.3]{ZF1}).
Since $\G$ is connected, ${\rm F}$ acts on $N_{\G}(1_0)$ faithfully.
By~\cite[Theorem~1.1 and Lemma~3.2]{ZF1}, we have the following proposition.

\begin{prop}\label{prop=biCay}
Let $\G=\BiCay(H, \emptyset, \emptyset, S)$ be a connected bi-Cayley graph over an abelian group $H$, and let $A=\Aut(\G)$. Then $N_{A}(R(H))=R(H)\lg{\rm F},\d_{\g,1,1}\rg$
with vertex stabilizer $(N_{A}(R(H)))_{1_0}={\rm F}$,
and $\G$ is isomorphic to the Cayley graph $\Cay(\Dih(H), \g S)$, where
$\Dih(H)=H\rtimes\lg \g\rg$.
\end{prop}

The following proposition is from~\cite[Theorem~1.1]{AHK}.

\begin{lem}\label{lem=cycliccover}
Let $n$ be a positive integer and $p$ a prime such that $p\geq5$. Let $\G$ be a connected pentavalent symmetric
bi-Cayley graph over $\mz_{np}$. Then $\G\cong\mathcal{CD}_{np}$,
as defined in Example~{\rm \ref{exam=n}}.
\end{lem}

Let $H=\lg x\rg\times\lg y\rg\times\lg z\rg=\mz_m\times\mz_{p^e}\times\mz_p$,
where $m$ and $e$ are two positive integers and $p$ is a prime such that $p\geq 5$ and $(m,p)=1$. In the remaining of this section, we always
let $\G=\BiCay(H,\emptyset, \emptyset, S)$ be a connected pentavalent bi-Cayley graph over $H$
such that $N_{{\rm Aut}(\G)}(R(H))$ is arc-transitive on $\G$.
Assume that $S=\{1,a,b,c,d\}$.
Then $H=\lg a,b,c,d\rg$. By Proposition~\ref{prop=biCay}, there exists a $\sigma_{\a,g}\in {\rm F}$ of order 5 permuting the neighborhood $\{1_1,a_1,b_1,c_1,d_1\}$ of $1_0$ in $\G$
cyclicly. One may assume that $1_1^{\sigma_{\a,g}}=a_1$, which implies that $g=a$ because $1_1^{\sigma_{\a,g}}=g_1$, and that $b_1=a_1^{\sigma_{\a,a}}$, $c_1=b_1^{\sigma_{\a,a}}$,
$d_1=c_1^{\sigma_{\a,a}}$ and  $1_1=d_1^{\sigma_{\a,a}}$. It follows that

\parbox{8cm}{
\begin{eqnarray*}
&& a^{\a}=ba^{-1},\ b^{\a}=ca^{-1},\ c^{\a}=da^{-1},\
d^{\a}=a^{-1}.
\end{eqnarray*}}\hfill
\parbox{1cm}{
\begin{eqnarray}
\label{eq=cover1-eq1}
\end{eqnarray}}

For $h\in H$, denote by $o(h)$ the order of $h$ in $H$.
Since $a^{\a}=ba^{-1}$ by Eq~(\ref{eq=cover1-eq1}),
$o(ba^{-1})=o(a^{\a})=o(a)$, forcing that $o(b)\di o(a)$.
Similarly, since $d^{\a}=a^{-1}$ and $c^{\a}=da^{-1}$,
we have $o(d)=o(a)$ and $o(c)\di o(a)$. Since $H=\lg a,b,c,d\rg$, we have $o(x)\di o(a)$ for any $x\in H$, and since $H=\mz_m\times\mz_{p^e}\times\mz_p$, we have $o(a)=mp^e$ and $|H:\langle a\rangle|=p$.

Suppose that $b\in\lg a\rg$, say $b=a^i$ for some integer $i$. By Eq~(\ref{eq=cover1-eq1}),
$a^{\a}=ba^{-1}=a^{i-1}\in \lg a\rg$ and $ca^{-1}=b^{\a}=(a^i)^{\a}=a^{i(i-1)}\in \lg a\rg$,
implying that $c\in\lg a\rg$. Similarly, $d\in \lg a\rg$ because $d=a\cdot c^{\a}$.
Since $H=\lg a,b,c,d\rg$, we have $H=\lg a\rg\cong\mz_{mp^e}$, a contradiction.
Hence, $b\notin \lg a\rg$, and since $|H:\lg a\rg|=p$, we have $H=\lg a,b\rg$ and $p\di o(b)$.

Let $A=\Aut(\G)$. Since $\G$ is $N_A(R(H))$-arc-transitive, ${\rm F}=N_A(R(H))_{1_0}$
acts transitively on $N_{\G}(1_0)$.
Let $\sigma_{\b,g}\in {\rm F}$ for some $\b\in\Aut(H)$ and $g\in H$ such that $1_1^{\sigma_{\b,g}}=1_1$. Then $1_1=(1^{\b}g)_1=g_1$,
forcing that $g=1$. Hence ${\rm F}_{1_1}=\{\sigma_{\b,1}~|~\b\in\Aut(H), S^{\b}=S\}$,
that is, ${\rm F}_{1_1}\cong\Aut(H,S)$.
By Proposition~\ref{prop=biCay}, $|N_A(R(H))|=2|H||{\rm F}|=2|H|\cdot|N_{\G}(1_0)||{\rm F}_{1_1}|=10|H||\Aut(H,S)|$.

\medskip

\f {\bf Observation:} $o(a)=mp^e$, $p\di o(b)$, $H=\lg a,b\rg$
and $|N_A(R(H))|=10|H||\Aut(H,S)|$.

\medskip

In the following two lemmas we consider the two cases: $e\geq 2$ and $e=1$, respectively.

\begin{lem}\label{lem=cover1}
If $e\geq2$, then $5\di (p-1)$, $\G\cong\mathcal{CGD}_{mp^e\times p}^{i}$ for some $1\leq i\leq 3$ and
$|N_A(R(H))|=10|H|$.
\end{lem}

\f {\bf Proof:} By Observation, $o(a)=mp^e$, $p\di o(b)$ and $H=\lg a,b\rg=\lg x,y,z\rg=\mz_m\times\mz_{p^e}\times\mz_p$, where $(m,p)=1$.
Then $H$ has an automorphism mapping $xy$ to $a$, and thus we
may assume $a=xy$, which implies that $b=x^{r+1}y^{\lambda+1}z^\iota$ for some $r+1\in\mz_m$, $\lambda+1\in\mz_{p^e}$ and $0\not=\iota\in\mz_p$ because $H=\lg a,b\rg$. Furthermore, $H$ has an automorphism fixing $x,y$ and mapping $z$ to $z^\iota$, and so we may assume $b=x^{r+1}y^{\lambda+1}z$. Let $c=x^iy^jz^{s}$ and $d=x^ky^{\ell}z^t$, where $i,k\in\mz_m$, $j,\ell\in\mz_{p^e}$ and $s,t\in\mz_p$.

Note that both $\lg x\rg=\mz_m$ and $\lg y,z\rg=\mz_{p^e}\times\mz_p$
are characteristic in $H$.
Since $a^{\a}=ba^{-1}$ by Eq~(\ref{eq=cover1-eq1}), that is,
$(xy)^{\a}=x^ry^{\lambda}z$, we have $x^{\a}=x^r$ and $y^{\a}=y^{\lambda}z$.
Since $(x^{r+1}y^{\lambda+1}z)^{\a}=b^{\a}=ca^{-1}=x^{i-1}y^{j-1}z^s$, we have
$z^{\a}=(x^{-r-1})^{\a}\cdot(y^{-\lambda-1})^{\a}\cdot(x^{i-1}y^{j-1}z^s)
=x^{-r^2-r-1+i}y^{-\lambda^2-\lambda-1+j}z^{s-\lambda-1}$, implying that
$z^{\a}=y^{-\lambda^2-\lambda-1+j}z^{s-\lambda-1}$ and

\parbox{8cm}{
\begin{eqnarray*}
&& -r^2-r-1+i=0~(\mod m),\\
&&-\lambda^2-\lambda-1+j=0~(\mod p^{e-1}).
\end{eqnarray*}}\hfill
\parbox{1cm}{
\begin{eqnarray}
\label{eq=cover1-eq3}\\ \label{eq=cover1-eq4}
\end{eqnarray}}

\f Similarly, since $c^{\a}=da^{-1}$ and $d^{\a}=a^{-1}$ by Eq~(\ref{eq=cover1-eq1}), we have
$x^{k-1}y^{\ell-1}z^t=c^{\a}=(x^iy^jz^{s})^{\a}=x^{ir}y^{\lambda j+s(-\lambda^2-\lambda-1+j)}z^{j+s(s-\lambda-1)}$ and
$x^{-1}y^{-1}=d^{\a}=(x^ky^{\ell}z^{t})^{\a}=x^{kr}y^{\lambda \ell+t(-\lambda^2-\lambda-1+j)}z^{\ell+t(s-\lambda-1)}$,
and by considering the powers of $x$, $y$ and $z$, we have the following Eqs~(\ref{eq=cover1-eq5})-(\ref{eq=cover1-eq10}).

\parbox{8cm}{
\begin{eqnarray*}
&& ir=k-1~(\mod m);\\
&& kr=-1~(\mod m);\\
&& \lambda j+s(-\lambda^2-\lambda-1+j)=\ell-1~(\mod p^e);\\
&& \lambda \ell+t(-\lambda^2-\lambda-1+j)=-1~(\mod p^e);\\
&& j+s(s-\lambda-1)=t~(\mod p);\\
&& \ell+t(s-\lambda-1)=0~(\mod p).
\end{eqnarray*}}\hfill
\parbox{1cm}{
\begin{eqnarray}
\label{eq=cover1-eq5}\\ \label{eq=cover1-eq6}\\ \label{eq=cover1-eq7}\\ \label{eq=cover1-eq8}
\\ \label{eq=cover1-eq9}\\ \label{eq=cover1-eq10}
\end{eqnarray}}

By Eq~(\ref{eq=cover1-eq3}), $i=r^2+r+1~(\mod m)$, and by
Eqs~(\ref{eq=cover1-eq5}) and
(\ref{eq=cover1-eq6}), $k=r^3+r^2+r+1~(\mod m)$ and
$r^4+r^3+r^2+r+1=0~(\mod m)$. It follows from Proposition~\ref{prop=x5=1}
that either $(r,m)\in \{(0,1),(1,5)\}$, or $r$ is an element in $\mz_m^*$ of order 5 and
$m=5^tp_1^{e_1}\cdots p_f^{e_f}$ with $t\leq 1$, $f\geq1$, $e_\iota\geq1$ and
$p_{\iota}$'s are distinct primes such that $5\di (p_\iota-1)$
for $1\leq \iota\leq f$.

Note that $e\geq2$. By Eq~(\ref{eq=cover1-eq4}),
$j=\lambda^2+\lambda+1~(\mod p^{e-1})$
and by Eqs~(\ref{eq=cover1-eq7})
and (\ref{eq=cover1-eq8}),  $\ell=\lambda^3+\lambda^2+\lambda+1~(\mod p^{e-1})$ and
$\lambda^4+\lambda^3+\lambda^2+\lambda+1=0~(\mod p^{e-1})$, implying $\lambda^5=1~(\mod p^{e-1})$. It follows from Proposition~\ref{prop=x5=1} that
either $(\lambda,p^{e-1})=(1,5)$, or $5\ |\ (p-1)$ and $\lambda$ is an element in $\mz_{p^{e-1}}^*$ of order 5, forcing that $\lambda\neq0$ and $\ell^{-1}=(-\lambda^4)^{-1}=-\lambda$. Furthermore, one may assume that $j=\lambda^2+\lambda+1+s_1p^{e-1}~(\mod p^e)$, $\ell=\lambda^3+\lambda^2+\lambda+1+s_2p^{e-1}~(\mod p^e)$ and $\lambda^4+\lambda^3+\lambda^2+\lambda+1=\iota p^{e-1}~(\mod p^e)$ for some $s_1,s_2,\iota\in\mz_p$.

In what follows all equations are considered in $\mz_p$,
unless otherwise stated. As $p\ |\ p^{e-1}$, the following equations are also true in $\mz_p$:
$$j=\lambda^2+\lambda+1,~\ell=\lambda^3+\lambda^2+\lambda+1,
~\lambda^4+\lambda^3+\lambda^2+\lambda+1=0,~\ell^{-1}=-\lambda.$$
By $s\times (\ref{eq=cover1-eq10})-t\times (\ref{eq=cover1-eq9})$,
$s=\ell^{-1}(jt-t^2)=-\lambda(jt-t^2)$, and by  Eq~(\ref{eq=cover1-eq10}),
we have $\lambda t^3-(\lambda^3+\lambda^2+\lambda)t^2-(\lambda+1)t+(\lambda^3+\lambda^2+\lambda+1)=0.$
Combined with $\lambda^4+\lambda^3+\lambda^2+\lambda+1=0$ and $\lambda\neq0$, we have
$(t-1)(t-\lambda)(t-\lambda^2)=0$,
which implies that $t=1$, $\lambda$ or $\lambda^2$. Recall that
$j=\lambda^2+\lambda+1$ and $s=-\lambda(jt-t^2)$. Thus
$(t,s)=(1,\lambda^4+\lambda+1)$, $(\lambda,\lambda^3+\lambda+1)$ or  $(\lambda^2,\lambda^2+\lambda+1)$.

Since $j=\lambda^2+\lambda+1+s_1p^{e-1}~(\mod p^e)$ and $\ell=\lambda^3+\lambda^2+\lambda+1+s_2p^{e-1}~(\mod p^e)$, by Eqs~(\ref{eq=cover1-eq7}) and (\ref{eq=cover1-eq8}) we have:

\parbox{8cm}{
\begin{eqnarray*}
\left\{
      \begin{array}{l}
         (\lambda+s)s_1p^{e-1}=s_2p^{e-1}~(\mod p^e)\\
         ts_1p^{e-1}+\lambda s_2p^{e-1}=-(\lambda^4+\lambda^3+\lambda^2+\lambda+1)~(\mod p^e)
      \end{array}
  \right.
\end{eqnarray*}}\hfill
\parbox{1cm}{
\begin{eqnarray}
\label{eq=cover1-eq11}
\end{eqnarray}}

\f Recall that either $(\lambda,p^{e-1})=(1,5)$ or $5\ |\ (p-1)$.
Suppose that $p^{e-1}=5$. Then $p=5$, $e=2$ and $(\lambda,s,t)=(1,3,1)$. By Eq~(\ref{eq=cover1-eq11}),
we have $5s_2=20s_1$ and $5^2s_1+5=0$ in $\mz_{5^2}$, a contradiction.
Hence $5\ |\ (p-1)$.
Again by Eq~(\ref{eq=cover1-eq11}), we have
$-(t+\lambda^2+\lambda s)s_1p^{e-1}=\iota p^{e-1}(\mod p^e)$,
where $\iota p^{e-1}=\lambda^4+\lambda^3+\lambda^2+\lambda+1$. Furthermore,

\parbox{8cm}{
\begin{eqnarray*}
\left\{
      \begin{array}{l}
         (t+\lambda^2+\lambda s)s_1=-\iota\\
          (t+\lambda^2+\lambda s)s_2=-\iota(\lambda+s)
              \end{array}
  \right.
\end{eqnarray*}}\hfill
\parbox{1cm}{
\begin{eqnarray}
\label{eq=cover1-eq12}
\end{eqnarray}}

Since $(t,s)=(1,\lambda^4+\lambda+1)$, $(\lambda,\lambda^3+\lambda+1)$ or  $(\lambda^2,\lambda^2+\lambda+1)$, we have
$t+\lambda^2+\lambda s=2\lambda^2+\lambda+2$, $\lambda^4+2\lambda^2+2\lambda$ or $\lambda^3+3\lambda^2+\lambda$,
respectively, and since $(2\lambda^2+\lambda+2)(\lambda^4+2\lambda^2+2\lambda)=6(\lambda^4+
\lambda^3+\lambda^2+\lambda)+1=-5$ and $(\lambda^3+3\lambda^2+\lambda)
(\lambda^4-2\lambda^3+\lambda^2)=\lambda^4+\lambda^3+\lambda^2+\lambda-4=-5$, we have
$(t+\lambda^2+\lambda s)^{-1}=-5^{-1}(\lambda^4+2\lambda^2+2\lambda)$, $-5^{-1}(2\lambda^2+\lambda+2)$
or $-5^{-1}(\lambda^4-2\lambda^3+\lambda^2)$, respectively. By Eq~(\ref{eq=cover1-eq12}),
$(s_1,s_2)=(5^{-1}\iota(\lambda^4+2\lambda^2+2\lambda), 5^{-1}\iota(-3\lambda^4+\lambda^3+2\lambda^2))$,
$(5^{-1}\iota(2\lambda^2+\lambda+2),5^{-1}\iota(-3\lambda^4+2\lambda^3+\lambda))$
or $(5^{-1}\iota(\lambda^4-2\lambda^3+\lambda^2),5^{-1}\iota(-2\lambda^4+\lambda^2+\lambda))$. It follows that
$$
\begin{array}{ll}
&a=xy,\ \ b=x^{r+1}y^{\lambda+1}z,\ \ (c,d)=\\
&(x^{r^2+r+1}y^{\lambda^2+\lambda+1+
5^{-1}(\lambda^4+2\lambda^2+2\lambda)\iota p^{e-1}}z^{\lambda^4+\lambda+1},
x^{r^3+r^2+r+1}y^{\lambda^3+\lambda^2+\lambda+1+
5^{-1}(-3\lambda^4+\lambda^3+2\lambda^2)\iota p^{e-1}}z),\\
&(x^{r^2+r+1}y^{\lambda^2+\lambda+1+5^{-1}(2\lambda^2+
\lambda+2)\iota p^{e-1}}z^{\lambda^3+\lambda+1},
x^{r^3+r^2+r+1}y^{\lambda^3+\lambda^2+\lambda+1+
5^{-1}(-3\lambda^4+2\lambda^3+\lambda)\iota p^{e-1}}z^{\lambda})\ {\rm or}\\
&(x^{r^2+r+1}y^{\lambda^2+\lambda+1+5^{-1}(\lambda^4-
2\lambda^3+\lambda^2)\iota p^{e-1}}z^{\lambda^2+\lambda+1},
x^{r^3+r^2+r+1}y^{\lambda^3+\lambda^2+\lambda+1+
5^{-1}(-2\lambda^4+\lambda^2+\lambda)\iota p^{e-1}}z^{\lambda^2}),
\end{array}$$
\f and $S=S_1$, $S_2$ or $S_3$, where

$$\begin{array}{llll}
&S_1=\{1,xy, x^{r+1}y^{\lambda+1}z, x^{r^2+r+1}y^{\lambda^2+\lambda+1+5^{-1}(\lambda^4+2\lambda^2+
2\lambda)\iota p^{e-1}}z^{\lambda^4+\lambda+1},\\
&~~~~~~~~x^{r^3+r^2+r+1}y^{\lambda^3+\lambda^2+\lambda+1+
5^{-1}(-3\lambda^4+\lambda^3+2\lambda^2)\iota p^{e-1}}z\},\\
&S_2=\{1,xy, x^{r+1}y^{\lambda+1}z, x^{r^2+r+1}y^{\lambda^2+\lambda+1+5^{-1}(2\lambda^2+
\lambda+2)\iota p^{e-1}}z^{\lambda^3+\lambda+1},\\
&~~~~~~~~x^{r^3+r^2+r+1}y^{\lambda^3+\lambda^2+\lambda+1+
5^{-1}(-3\lambda^4+2\lambda^3+\lambda)\iota p^{e-1}}z^{\lambda}\},\\
&S_3=\{1,xy, x^{r+1}y^{\lambda+1}z, x^{r^2+r+1}y^{\lambda^2+\lambda+1+5^{-1}(\lambda^4-2\lambda^3+
\lambda^2)\iota p^{e-1}}z^{\lambda^2+\lambda+1},\\
&~~~~~~~~x^{r^3+r^2+r+1}y^{\lambda^3+\lambda^2+\lambda+1+
5^{-1}(-2\lambda^4+\lambda^2+\lambda)\iota p^{e-1}}z^{\lambda^2}\}.
\end{array}$$

Since $x^5=1~(\mod p^e)$ implies that  $x^5=1~(\mod p^{e-1})$, there exists $f\in\mz_p$ such that $\lambda_1=\lambda+fp^{e-1}$ is an element of order $5$ in $\mz_{p^e}^*$.
Then $\lambda=\lambda_1-fp^{e-1}$, $\lambda_1^5=1$ and $\lambda_1^4+\lambda_1^3+\lambda_1^2+\lambda_1+1=0$ in $\mz_{p^e}$. Hence $\iota p^{e-1}=\lambda^4+\lambda^3+\lambda^2+\lambda+1=
(\lambda_1-fp^{e-1})^4+(\lambda_1-fp^{e-1})^3+(\lambda_1-fp^{e-1})^2+(\lambda_1-fp^{e-1})+1
=-(4\lambda_1^3+3\lambda_1^2+2\lambda_1+1)fp^{e-1}$
in $\mz_{p^{e}}$, and thus
$$\begin{array}{llll}
&S_1=\{1,xy, x^{r+1}y^{\lambda+1}z, x^{r^2+r+1}y^{\lambda^2+\lambda+1+5^{-1}(\lambda^4+2\lambda^2+
2\lambda)\iota p^{e-1}}z^{\lambda^4+\lambda+1},\\
&~~~~~~~~x^{r^3+r^2+r+1}y^{\lambda^3+\lambda^2+\lambda+1+
5^{-1}(-3\lambda^4+\lambda^3+2\lambda^2)\iota p^{e-1}}z\},\\
&~~~=\{1,xy, x^{r+1}y^{\lambda_1+1}y^{-fp^{e-1}}z, x^{r^2+r+1}y^{\lambda_1^2+\lambda_1+1}y^{-(\lambda_1^4+\lambda_1+1)fp^{e-1}}z^{\lambda_1^4+\lambda_1+1},\\
&~~~~~~~~x^{r^3+r^2+r+1}y^{\lambda_1^3+\lambda_1^2+\lambda_1+1}y^{-fp^{e-1}}z\}.
\end{array}$$

\f Let $\varphi$ be the automorphism of $H$ induced by $x\mapsto x$, $y\mapsto y$ and $z\mapsto y^{fp^{e-1}}z$. Then $(S_1)^{\varphi}=\{1,xy, x^{r+1}y^{\lambda_1+1}z, x^{r^2+r+1}y^{\lambda_1^2+\lambda_1+1}z^{\lambda_1^4+\lambda_1+1},
x^{r^3+r^2+r+1}y^{\lambda_1^3+\lambda_1^2+\lambda_1+1}z\}$.
Since $\BiCay(H,\emptyset,\emptyset, S_1)\cong \BiCay(H,\emptyset,\emptyset, S_1^{\varphi})$,
we may assume that $\lambda=\lambda_1$ is an element of order $5$ in $\mz_{p^e}^*$, and
$$S_1=\{1,xy, x^{r+1}y^{\lambda+1}z, x^{r^2+r+1}y^{\lambda^2+\lambda+1}z^{\lambda^4+\lambda+1},
x^{r^3+r^2+r+1}y^{\lambda^3+\lambda^2+\lambda+1}z\}.$$
Similarly, we can also assume that
$$\begin{array}{llll}
&S_2=\{1,xy, x^{r+1}y^{\lambda+1}z, x^{r^2+r+1}
y^{\lambda^2+\lambda+1}z^{\lambda^3+\lambda+1},
x^{r^3+r^2+r+1}y^{\lambda^3+\lambda^2+\lambda+1}z^{\lambda}\},\\
&S_3=\{1,xy, x^{r+1}y^{\lambda+1}z,
x^{r^2+r+1}y^{\lambda^2+\lambda+1}z^{\lambda^2+\lambda+1},
x^{r^3+r^2+r+1}y^{\lambda^3+\lambda^2+\lambda+1}z^{\lambda^2}\}.
\end{array}$$
\f By Proposition~\ref{prop=biCay} and Example~\ref{exam=n*p}, $\G=\BiCay(H,\emptyset,\emptyset, S)\cong\mathcal{CGD}_{mp^e\times p}^{i}$ with $1\leq i\leq3$.

Note that $|N_A(R(H))|=10|H||\Aut(H,S)|$ by Observation.
For $S_1$, let $\b\in \Aut(H,S_1)$. Then $S_1^{\b}=S_1$.
Since $\lg y,z\rg$ is characteristic in $H=\lg x\rg\times\lg y\rg\times\lg z\rg=\mz_m\times\mz_{p^e}\times\mz_p$, we have
$$\{y, y^{\lambda+1}z, y^{\lambda^2+\lambda+1}z^{\lambda^4+\lambda+1},
y^{\lambda^3+\lambda^2+\lambda+1}z\}^{\b}=\{y, y^{\lambda+1}z, y^{\lambda^2+\lambda+1}z^{\lambda^4+\lambda+1},
y^{\lambda^3+\lambda^2+\lambda+1}z\}.$$
It follows that $y^{\b}=y^sz^t$ with $(s,t)=(1,0)$,
$(\lambda+1,1)$, $(\lambda^2+\lambda+1,\lambda^4+\lambda+1)$,
or $(\lambda^3+\lambda^2+\lambda+1,1)$.
Furthermore, we have $(y\cdot y^{\lambda+1}z\cdot y^{\lambda^2+\lambda+1}z^{\lambda^4+\lambda+1}\cdot
y^{\lambda^3+\lambda^2+\lambda+1}z)^{\b}=y\cdot y^{\lambda+1}z\cdot y^{\lambda^2+\lambda+1}z^{\lambda^4+\lambda+1}\cdot
y^{\lambda^3+\lambda^2+\lambda+1}z,$
that is, $(y^{\b}y^{-1})^{\lambda^3+2\lambda^2+3\lambda+4}=(z^{-\lambda^4-\lambda-3})^{\b}z^{\lambda^4+\lambda+3}$.
In particular, $(y^{\b}y^{-1})^{(\lambda^3+2\lambda^2+3\lambda+4)p}=1$.
Note that $\lambda^4+\lambda^3+\lambda^2+\lambda+1=0$ in $\mz_{p^e}$
implies that $\lambda^4+\lambda^3+\lambda^2+\lambda+1=0$ in $\mz_{p}$.
If $\lambda^3+2\lambda^2+3\lambda+4=0$ in $\mz_p$,
then $\lambda^3=-2\lambda^2-3\lambda-4$,
$\lambda^4=\lambda\cdot \lambda^3=\lambda^2+2\lambda+8$,
and thus $0=\lambda^4+\lambda^3+\lambda^2+\lambda+1=5$, contrary to
$5\di (p-1)$. Hence $\lambda^3+2\lambda^2+3\lambda+4\neq 0$ in $\mz_p$
and $(y^{\b}y^{-1})^{p}=1$.

Suppose that $(s,t)\neq (1,0)$. Then $y^\b y^{-1}=y^{s-1}z^t$ with $s-1=\lambda,\lambda^2+\lambda$ or $\lambda^3+\lambda^2+\lambda$. Since $\lambda^4+\lambda^3+\lambda^2+\lambda+1=0$ in $\mz_p$, we have $\lambda\neq 0,-1$ and thus $(s-1,p)=1$. This implies that $y^\b y^{-1}=y^{s-1}z^t$ has order $p^e$, and since $e\geq 2$, we have  $(y^\b y^{-1})^p\not=1$, a contradiction. Hence $(s,t)=(1,0)$, that is, $y^\b=y$. It follows that $(y^{\lambda+1}z)^{\b}=y^{\lambda+1}z^{\b}\in\{y^{\lambda+1}z, y^{\lambda^2+\lambda+1}z^{\lambda^4+\lambda+1},
y^{\lambda^3+\lambda^2+\lambda+1}z\}$,
and thus $z^{\b}\in\{z,y^{\lambda^2}z^{\lambda^4+\lambda+1},
y^{\lambda^3+\lambda^2}z\}$. If $z^{\b}=y^{\lambda^2}z^{\lambda^4+\lambda+1}$
or $y^{\lambda^3+\lambda^2}z$, then $(y^{\lambda^2})^p=1$ or
$(y^{\lambda^3+\lambda^2})^p=1$. It forces that
$\lambda^2=0$ or $\lambda^3+\lambda^2=0$ in $\mz_{p^{e-1}}$,
and $\lambda=0,-1$, a contradiction. Hence $z^{\b}=z$.
Noting that $\lg x\rg$ is characteristic in $H$,
we have $(xy)^{\b}=x^{\b}y\in S_1^{\b}=S_1$. Then it is easy to check that $(xy)^{\b}=xy$
and thus $x^{\b}=x$.
It implies that $\b$ is the identity automorphism. Hence $|\Aut(H,S_1)|=1$
and $|N_A(R(H))|=10|H|$. By a similar argument as above, for $S_2$ and $S_3$, we also have $|\Aut(H,S_2)|=|\Aut(H,S_3)|=1$ and $|N_A(R(H))|=10|H|$.  \hfill$\blacksquare$

\begin{lem}\label{lem=cover2}
If $e=1$, that is, $H\cong\mz_m\times\mz_p\times\mz_p$, then one of the following holds:

\begin{enumerate}

\item [{\rm (1)}] $p=5$ or $5\di (p\pm1)$ and
$\G\cong\mathcal{CGD}_{mp\times p}^{4}$ as defined in
Example~{\rm \ref{exam=n*p*p-1}}. Furthermore,
$|N_A(R(H))|=10|H|$ if $m\neq 1,5$;
$|N_A(R(H))|=20|H|$ if $m=5$;
$|N_A(R(H))|=20|H|$ if $m=1$ and $p\neq5$;
and $|N_A(R(H))|=40|H|$ if $m=1$ and $p=5$;

\item [{\rm (2)}] $5\di (p-1)$, $\G\cong\mathcal{CGD}_{mp\times p}^{5}$
as defined in Example~{\rm\ref{exam=n*p*p-2}}
and $|N_A(R(H))|=10|H|$.
\end{enumerate}
\end{lem}

\f {\bf Proof:} Note that $(m,p)=1$. By Observation, we have $o(a)=mp$, $p\di o(b)$ and $H=\lg a,b\rg=\lg x,y,z\rg=\mz_m\times\mz_{p}\times\mz_p$.
Then $H$ has an automorphism mapping $xy$ to $a$, and
we may assume $a=xy$, implying that $b=x^{r+1}y^{\lambda}z^\iota$ for some $r+1\in\mz_m$, $\iota,\lambda\in\mz_{p}$ and $\iota\neq 0$ because $H=\lg a,b\rg$.
The group $H$ also has an automorphism fixing $x,y$ and mapping $z$ to $y^{\lambda}z^\iota$, and we may further assume $b=x^{r+1}z$. Let
$c=x^iy^jz^{s}$ and $d=x^ky^{\ell}z^t$, where $i,k\in\mz_m$, $j,\ell,s,t\in\mz_{p}$.

By Eq~(\ref{eq=cover1-eq1}), $a^{\a}=ba^{-1}$, that is, $(xy)^{\a}=x^{r}y^{-1}z$.
Since both $\lg x\rg$ and $\lg y,z\rg$ are characteristic in $H$, we have
$x^{\a}=x^r$ and $y^{\a}=y^{-1}z$.
Again by Eq~(\ref{eq=cover1-eq1}), since
$(x^{r+1}z)^{\a}=b^{\a}=ca^{-1}=x^{i-1}y^{j-1}z^s$,
we have $z^{\a}=(x^{-r-1})^{\a}\cdot b^{\a}=x^{-r^2-r-1+i}y^{j-1}z^s$,
implying that $z^{\a}=y^{j-1}z^s$ and

\parbox{8cm}{
\begin{eqnarray*}
&& -r^2-r-1+i=0~(\mod m).
\end{eqnarray*}}\hfill
\parbox{1cm}{
\begin{eqnarray}
\label{eq=cover2-eq1}
\end{eqnarray}}

\f Note that
$x^{k-1}y^{\ell-1}z^t=da^{-1}=c^{\a}=
(x^iy^jz^s)^{\a}=(x^r)^i(y^{-1}z)^j(y^{j-1}z^s)^s=x^{ri}y^{-j+s(j-1)}z^{j+s^2}$
and $x^{-1}y^{-1}=a^{-1}=d^{\a}=(x^ky^{\ell}z^t)^{\a}=
(x^r)^k(y^{-1}z)^{\ell}(y^{j-1}z^s)^t=x^{rk}y^{-\ell+(j-1)t}z^{st+\ell}$.
Considering the powers of $x$, $y$ and $z$, we have
Eqs~(\ref{eq=cover2-eq2})-(\ref{eq=cover2-eq7}).
As shown in these equations, in what follows all equations are considered in $\mz_p$, unless otherwise stated.

\parbox{8cm}{
\begin{eqnarray*}
&& k-1=ri~(\mod m);\\
&& \ell-1=-j+s(j-1);\\
&& t=j+s^2;\\
&& -1=rk~(\mod m);\\
&& -1=-\ell+(j-1)t;\\
&& 0=st+\ell.
\end{eqnarray*}}\hfill
\parbox{1cm}{
\begin{eqnarray}
\label{eq=cover2-eq2}\\ \label{eq=cover2-eq3}\\ \label{eq=cover2-eq4}\\ \label{eq=cover2-eq5}\\
\label{eq=cover2-eq6}\\ \label{eq=cover2-eq7}
\end{eqnarray}}

By Eq~(\ref{eq=cover2-eq1}),
we have $i=r^2+r+1~(\mod m)$ and by Eqs~(\ref{eq=cover2-eq2}) and
(\ref{eq=cover2-eq5}), $k=r^3+r^2+r+1~(\mod m)$ and $r^4+r^3+r^2+r+1=0~(\mod m)$.
It follows from Proposition~\ref{prop=x5=1} that either $(r,m)\in \{(0,1),(1,5)\}$ or $r$ is an element of order $5$ in $\mz_m^*$ and the prime decomposition of $m$ is $5^tp_1^{e_1}\cdots p_f^{e_f}$ with
$t\leq1$, $f\geq1$, $e_{\iota}\geq 1$ and $5\di (p_{\iota}-1)$ for $1\leq \iota\leq f$.

By Eq~(\ref{eq=cover2-eq4}), $t=j+s^2$,
and by Eqs~(\ref{eq=cover2-eq3}), (\ref{eq=cover2-eq6}) and (\ref{eq=cover2-eq7}), $\ell=1-j+s(j-1)$,
$\ell=1+(j-1)t=1+(j-1)(j+s^2)$ and $\ell=-st=-sj-s^3$. It follows

\parbox{8cm}{
\begin{eqnarray*}
&& j^2+(s^2-s)j-(s^2-s)=0;\\
&& (2s-1)j+s^3-s+1=0.
\end{eqnarray*}}\hfill
\parbox{1cm}{
\begin{eqnarray}
\label{eq=cover2-eq8}\\ \label{eq=cover2-eq9}
\end{eqnarray}}

\f By Eq~(\ref{eq=cover2-eq8}), $(2s-1)^2j^2+(2s-1)^2(s^2-s)j-(2s-1)^2(s^2-s)=0$, and
since $(2s-1)j=-(s^3-s+1)$, we have $s^6-3s^5+5s^4-5s^3+2s-1=0$, that is,
$(s^2-s-1)(s^4-2s^3+4s^2-3s+1)=0$. Hence, either $s^2-s-1=0$ or $s^4-2s^3+4s^2-3s+1=0$.

\medskip

\f {\bf Case~1:} $s^2-s-1=0$.

Let $\lambda=2s-1$. Then $s=2^{-1}(1+\lambda)$ and $\lambda^2=5$,
and thus $(\lambda, p)=(0,5)$ or $5\ |\ (p\pm 1)$ by~\cite[Example~4.6]{YF}.
By Eqs~(\ref{eq=cover2-eq8}) and (\ref{eq=cover2-eq9}), $j^2+j-1=0$ and $(2s-1)j+(s+2)=0$.

For $(\lambda, p)=(0,5)$, $j^2+j-1=0$ implies that $j=2=-2^{-1}(1+\lambda)$.
For $5\di (p\pm 1)$, we have $\lambda\neq 0$, and since
$2s-1=\lambda$ and $(2s-1)j+(s+2)=0$,
we have $j=-(2s-1)^{-1}(s+2)=-\lambda^{-1}\cdot 2^{-1}(\lambda+5)=-2^{-1}(1+\lambda)$
(note that $5=\lambda^2$).
It follows from Eqs~(\ref{eq=cover2-eq4}) and~(\ref{eq=cover2-eq7}) that
$t=j+s^2=1$  and $\ell=-st=-2^{-1}(1+\lambda)$.
Recall that $i=r^2+r+1~(\mod m)$ and $k=i^3+i^2+i+1~(\mod m)$. Hence
$c=x^{r^2+r+1}y^{-2^{-1}(1+\lambda)}z^{2^{-1}(1+\lambda)}$
and $d=x^{r^3+r^2+r+1}y^{-2^{-1}(1+\lambda)}z$. Now,
$$S=\{1,xy,x^{r+1}z,x^{r^2+r+1}y^{-2^{-1}(1+\lambda)}
z^{2^{-1}(1+\lambda)},x^{r^3+r^2+r+1}y^{-2^{-1}(1+\lambda)}z\}.$$
By Proposition~\ref{prop=biCay} and Example~\ref{exam=n*p*p-1},
$\G=\BiCay(H,\emptyset,\emptyset,S)\cong\mathcal{CGD}_{mp\times p}^{4}$.

For $(m,p)=(1,5)$, we have $\lambda=0$ and $S=\{1,y,z,y^{-3}z^{3},y^{-3}z\}$.
By MAGMA~\cite{BCP}, $|N_A(R(H))|=40|H|$.
Assume that $(m,p)\neq(1,5)$. To determine the order of $N_A(R(H))$, by Observation, we only need to determine $|\Aut(H,S)|$. Let $\b\in \Aut(H,S)$. Then

\parbox{8cm}{
\begin{eqnarray*}
&&\{1,xy, x^{r+1}z, x^{r^2+r+1}y^{-2^{-1}(1+\lambda)}
z^{2^{-1}(1+\lambda)}, x^{r^3+r^2+r+1}y^{-2^{-1}(1+\lambda)}z\}^{\b}\\
&&=\{1,xy, x^{r+1}z, x^{r^2+r+1}y^{-2^{-1}(1+\lambda)}
z^{2^{-1}(1+\lambda)}, x^{r^3+r^2+r+1}y^{-2^{-1}(1+\lambda)}z\}.
\end{eqnarray*}}\hfill
\parbox{1cm}{
\begin{eqnarray}
\label{eq=cover2-eqS}
\end{eqnarray}}

Let $m\neq 1$ or 5. Then $r\neq 0,\pm1$.
Since $\lg x\rg$ is a characteristic subgroup of $H$,
we have $\{x,x^{r+1},x^{r^2+r+1},x^{r^3+r^2+r+1}\}^{\b}=\{x,x^{r+1},x^{r^2+r+1},x^{r^3+r^2+r+1}\}$
by Eq~(\ref{eq=cover2-eqS}), and thus $(x\cdot x^{r+1}\cdot x^{r^2+r+1}\cdot x^{r^3+r^2+r+1})^{\b}
=x\cdot x^{r+1}\cdot x^{r^2+r+1}\cdot x^{r^3+r^2+r+1}$, that is, $(x^{\b}x^{-1})^{r^3+2r^2+3r+4}=1$.
Suppose that $x^{\b}=x^{r+1}$, $x^{r^2+r+1}$ or $x^{r^3+r^2+r+1}$.
Then $1=(x^{\b}x^{-1})^{r^3+2r^2+3r+4}=x^{r(r^3+2r^2+3r+4)}$,
$x^{r(r+1)(r^3+2r^2+3r+4)}$ or $x^{r(r^2+r+1)(r^3+2r^2+3r+4)}$,
implying that $r(r^3+2r^2+3r+4)=0$, $r(r+1)(r^3+2r^2+3r+4)=0$
or $r(r^2+r+1)(r^3+2r^2+3r+4)=0$ in $\mz_m$.
Note that $r\in\mz_m^*$.
Suppose that $r+1=0~(\mod q)$ for some prime divisor $q$ of $m$.
Since $r^4+r^3+r^2+r+1=0~(\mod m)$,
we have $r^4+r^3+r^2+r+1=0~(\mod q)$, and since $r=-1~(\mod q)$,
we have $1=0~(\mod q)$, a contradiction.
Hence $r+1\in\mz_m^*$, and $r^2+r+1=-r^3(r+1)\in\mz_m^*$. Thus $r^3+2r^2+3r+4=0$ in $\mz_m$.
Then $r^3=-2r^2-3r-4$, $r^4=r\cdot r^3=r^2+2r+8$,
and $0=r^4+r^3+r^2+r+1=5$ in $\mz_m$, forcing that $m=1$ or 5, a contradiction.
Hence $x^{\b}=x$. Since $r\neq0,\pm1$, any two elements in $\{1, r+1,r^2+r+1, r^3+r^2+r+1\}$ are not equal. It follows from Eq~(\ref{eq=cover2-eqS}) that $(xy)^{\b}=xy$ and $(x^{r+1}z)^{\b}=x^{r+1}z$, forcing that $y^{\b}=y$ and $z^{\b}=z$.
It implies $|\Aut(H,S)|=1$, and thus $|N_A(R(H))|=10|H||\Aut(H,S)|=10|H|$.

Let $m=1$ or 5. Then $r=0$ or 1, and $p\neq5$
because $(m,p)\neq (1,5)$ and $(m,p)=1$.
Since $\lambda^2=5$, we have
$\lambda\neq0,\pm1$ and $2^{-1}(5+\lambda)\cdot 10^{-1}(5-\lambda)=1$.
By Eq~(\ref{eq=cover2-eqS}), we have
$\{y,z,y^{-2^{-1}(1+\lambda)}
z^{2^{-1}(1+\lambda)},y^{-2^{-1}(1+\lambda)}z\}^{\b}=\{y,z,y^{-2^{-1}(1+\lambda)}
z^{2^{-1}(1+\lambda)},y^{-2^{-1}(1+\lambda)}z\}$ because $\lg y,z\rg$
is characteristic in $H$.
Hence
$(y\cdot z\cdot y^{-2^{-1}(1+\lambda)}z^{2^{-1}(1+\lambda)}\cdot
 y^{-2^{-1}(1+\lambda)}z)^{\b}=y\cdot z\cdot y^{-2^{-1}(1+\lambda)}z^{2^{-1}(1+\lambda)}\cdot y^{-2^{-1}(1+\lambda)}z$,
that is, $(y^{-\lambda}z^{2^{-1}(5+\lambda)})^{\b}=
y^{-\lambda}z^{2^{-1}(5+\lambda)}$,
and so $z^{\b}=
[(y^{\lambda})^{\b}\cdot y^{-\lambda}z^{2^{-1}(5+\lambda)}]^{10^{-1}(5-\lambda)}
=(y^{\b})^{2^{-1}(\lambda-1)}y^{2^{-1}(1-\lambda)}z\in\{y,z,y^{-2^{-1}(1+\lambda)}
z^{2^{-1}(1+\lambda)},y^{-2^{-1}(1+\lambda)}z\}$.

Suppose $(xy)^{\b}=x^{r+1}z$. Then $x^{\b}=x^{r+1}$, $y^{\b}=z$,
and $z^{\b}=(y^{\b})^{2^{-1}(\lambda-1)}y^{2^{-1}(1-\lambda)}z
=y^{2^{-1}(1-\lambda)}z^{2^{-1}(1+\lambda)}$. Since $z^{\b}\in \{y,z,y^{-2^{-1}(1+\lambda)}
z^{2^{-1}(1+\lambda)},y^{-2^{-1}(1+\lambda)}z\}$,
we have $2^{-1}(1-\lambda)=1$, 0, or $-2^{-1}(1+\lambda)$
by considering the power of $y$. It forces that $\lambda=\pm1$ or $2=0$, a contradiction.

Suppose $(xy)^{\b}=x^{r^2+r+1}y^{-2^{-1}(1+\lambda)}z^{2^{-1}(1+\lambda)}$.
Then $x^{\b}=x^{r^2+r+1}$, $y^{\b}=y^{-2^{-1}(1+\lambda)}z^{2^{-1}(1+\lambda)}$,
and $z^{\b}=(y^{-2^{-1}(1+\lambda)}z^{2^{-1}(1+\lambda)})^{2^{-1}(\lambda-1)}y^{2^{-1}(1-\lambda)}z=
y^{-2^{-1}(1+\lambda)}z^2$. It is easy to check that $z^{\b}\notin\{y,z,y^{-2^{-1}(1+\lambda)}
z^{2^{-1}(1+\lambda)},y^{-2^{-1}(1+\lambda)}z\}$ by considering the power of $z$,
a contradiction.

By Eq~(\ref{eq=cover2-eqS}), either $(xy)^{\b}=xy$ or $x^{r^3+r^2+r+1}y^{-2^{-1}(1+\lambda)}z$.
For the former, $x^{\b}=x$, $y^{\b}=y$ and $z^{\b}=y^{2^{-1}(\lambda-1)}y^{2^{-1}(1-\lambda)}z=z$, implying that $\b$ is the identity automorphism.
For the latter,
we have $x^{\b}=x^{r^3+r^2+r+1}$, $y^{\b}=y^{-2^{-1}(1+\lambda)}z$,
and $z^{\b}=(y^{-2^{-1}(1+\lambda)}z)^{2^{-1}(\lambda-1)}y^{2^{-1}(1-\lambda)}z
=y^{-2^{-1}(\lambda+1)}z^{2^{-1}(\lambda+1)}$.
Noting that $r=0$ or 1, we have $r^3+r^2+r+1=1$ or $-1$,
and thus $\b$ has order 2. Hence $|\Aut(G,S)|=2$, and by Observation, $|N_A(R(H))|=10|H||\Aut(H,S)|=20|H|$.

\medskip

\f {\bf Case~2:} $s^4-2s^3+4s^2-3s+1=0$.

By Case 1, we may assume that $s^2-s-1\not=0$. If $p=5$, then
$s^4-2s^3+4s^2-3s+1=0$ implies that $s=3$ and thus $s^2-s-1=0$,
a contradiction. Hence $p\neq5$. By~\cite[Lemma~5.4, Case 2]{YF},
we have $5\ |\ (p-1)$ and $s=2^{-1}(1+\lambda)$,
where $\lambda^4+10\lambda^2+5=0$ and $\lambda\neq0,\pm 1$.

Since $s^4-2s^3+4s^2-3s+1=0$, we have
$(2s-1)(8s^3-12s^2+26s-11)=-5$, and since $p\neq5$,
we have $(2s-1)^{-1}=-5^{-1}(8s^3-12s^2+26s-11)$.
Noting that $s^4=2s^3-4s^2+3s-1$, we have $s^5=-5s^2+5s-2$ and $s^6=-5s^3+5s^2-2s$.
By Eq~(\ref{eq=cover2-eq9}), $j=-(2s-1)^{-1}(s^3-s+1)=5^{-1}(8s^3-12s^2+26s-11)(s^3-s+1)=
s^3-2s^2+3s-1=8^{-1}(\lambda^3-\lambda^2+7\lambda+1)$ and
by Eqs~(\ref{eq=cover2-eq4}) and (\ref{eq=cover2-eq7}), $t=j+s^2=s^3-s^2+3s-1=8^{-1}(\lambda^3+\lambda^2+11\lambda+3)$ and
$\ell=-st=-s^3+s^2-2s+1=-8^{-1}(\lambda^3+\lambda^2+7\lambda-1)$.
It follows that
$$S=\{1,xy,x^{r+1}z,x^{r^2+r+1}y^{8^{-1}(\lambda^3-\lambda^2+7\lambda+1)}
z^{2^{-1}(1+\lambda)},x^{r^3+r^2+r+1}y^{-8^{-1}(\lambda^3+\lambda^2+7\lambda-1)}
z^{8^{-1}(\lambda^3+\lambda^2+11\lambda+3)}\}.
$$
By Proposition~\ref{prop=biCay} and Example~\ref{exam=n*p*p-2},
$\G=\BiCay(H,\emptyset,\emptyset,S)\cong\mathcal{CGD}_{mp\times p}^{5}$.

Let $\b\in\Aut(H,S)$. Then $S^{\b}=S$. Since $\lg x\rg$ and $\lg y,z\rg$ are characteristic in $H$, we have

\parbox{8cm}{
\begin{eqnarray*}
&&\{y,z,y^{8^{-1}(\lambda^3-\lambda^2+7\lambda+1)}
z^{2^{-1}(1+\lambda)},y^{-8^{-1}(\lambda^3+\lambda^2+7\lambda-1)}
z^{8^{-1}(\lambda^3+\lambda^2+11\lambda+3)}\}^{\b}\\
&&=\{y,z,y^{8^{-1}(\lambda^3-\lambda^2+7\lambda+1)}
z^{2^{-1}(1+\lambda)},y^{-8^{-1}(\lambda^3+\lambda^2+7\lambda-1)}
z^{8^{-1}(\lambda^3+\lambda^2+11\lambda+3)}\}.
\end{eqnarray*}}\hfill
\parbox{1cm}{
\begin{eqnarray}
\label{eq=cover2-eqS2}
\end{eqnarray}}

\f It implies that
$(y\cdot z\cdot y^{8^{-1}(\lambda^3-\lambda^2+7\lambda+1)}
z^{2^{-1}(1+\lambda)}\cdot y^{-8^{-1}(\lambda^3+\lambda^2+7\lambda-1)}
z^{8^{-1}(\lambda^3+\lambda^2+11\lambda+3)})^{\b}=y\cdot z\cdot y^{8^{-1}(\lambda^3-\lambda^2+7\lambda+1)}
z^{2^{-1}(1+\lambda)}\cdot y^{-8^{-1}(\lambda^3+\lambda^2+7\lambda-1)}
z^{8^{-1}(\lambda^3+\lambda^2+11\lambda+3)}$,
that is, $(y^{10-2\lambda^2}z^{\lambda^3+\lambda^2+15\lambda+15})^{\b}=
y^{10-2\lambda^2}z^{\lambda^3+\lambda^2+15\lambda+15}$.
Set $f=-320^{-1}(\lambda^3-\lambda^2+15\lambda-15)$.
Recall that $\lambda^4=-10\lambda^2-5$,
$\lambda^5=-10\lambda^3-5\lambda$ and
$\lambda^6=95\lambda^2+50$.
Then $f\cdot(\lambda^3+\lambda^2+15\lambda+15)=1$
and
$[(y^{10-2\lambda^2}z^{\lambda^3+\lambda^2+15\lambda+15})^{\b}
]^{f}=(y^{\b})^{2^{-1}(1-\lambda)}z^{\b}=
(y^{10-2\lambda^2}z^{\lambda^3+\lambda^2+15\lambda+15})^{f}=y^{2^{-1}(1-\lambda)}z$,
yielding that $z^{\b}=(y^{\b}y^{-1})^{2^{-1}(\lambda-1)}z$.
Note that $5\di (p-1)$ and $\lambda^4+10\lambda^2+5=0$
imply that $\lambda\neq 0,\pm1$ and $\lambda^2\neq\pm 1,-3$.

Suppose $y^{\b}=z$. Then $z^{\b}=(zy^{-1})^{2^{-1}(\lambda-1)}z=y^{2^{-1}(1-\lambda)}z^{2^{-1}(\lambda+1)}$.
Since $\lambda\neq -1$, we have $2^{-1}(\lambda+1)\neq0$,
forcing that $z^{\b}=y^{2^{-1}(1-\lambda)}z^{2^{-1}(\lambda+1)}\neq y$. By Eq~(\ref{eq=cover2-eqS2}),
$z^{\b}=y^{8^{-1}(\lambda^3-\lambda^2+7\lambda+1)}
z^{2^{-1}(1+\lambda)}$ or $y^{-8^{-1}(\lambda^3+\lambda^2+7\lambda-1)}
z^{8^{-1}(\lambda^3+\lambda^2+11\lambda+3)}$,
and by considering the power of $y$, we have $2^{-1}(1-\lambda)=8^{-1}(\lambda^3-\lambda^2+7\lambda+1)$
or $-8^{-1}(\lambda^3+\lambda^2+7\lambda-1)$, that is
$\lambda^3=\lambda^2-11\lambda+3$ or $\lambda^3=-\lambda^2-3\lambda-3$.
Then $\lambda^4=\lambda\cdot\lambda^3=-10\lambda^2-8\lambda+3$
or $-2\lambda^2+3$, and $0=\lambda^4+10\lambda^2+5=8-8\lambda$ or $8\lambda^2+8$, yielding that $\lambda=1$
or $\lambda^2=-1$, a contradiction.

Suppose $y^{\b}=y^{8^{-1}(\lambda^3-\lambda^2+7\lambda+1)}
z^{2^{-1}(1+\lambda)}$. Then $z^{\b}=y^{-8^{-1}(\lambda^3+\lambda^2+7\lambda-1)}z^{4^{-1}(\lambda^2+3)}$.
Since $\lambda^2\neq-3$ or 1,
we have $4^{-1}(\lambda^2+3)\neq 0$ or 1, forcing that
$z^{\b}\neq y$ or $z$. It follows from Eq~(\ref{eq=cover2-eqS2})
that $z^{\b}=y^{-8^{-1}(\lambda^3+\lambda^2+7\lambda-1)}
z^{8^{-1}(\lambda^3+\lambda^2+11\lambda+3)}$
and $4^{-1}(\lambda^2+3)=8^{-1}(\lambda^3+\lambda^2+11\lambda+3)$
by considering the power of $z$,
yielding that $\lambda^3=\lambda^2-11\lambda+3$.
Thus $\lambda^4=\lambda\cdot\lambda^3=-10\lambda^2-8\lambda+3$ and
$0=\lambda^4+10\lambda^2+5=-8\lambda+8$, contrary to that $\lambda\neq-1$.

Suppose $y^{\b}=y^{-8^{-1}(\lambda^3+\lambda^2+7\lambda-1)}
z^{8^{-1}(\lambda^3+\lambda^2+11\lambda+3)}$. Then $z^{\b}=y^{4^{-1}(\lambda^3+3)}z^{2^{-1}(1-\lambda)}$,
and since $\lambda\neq\pm1,0$, we have $z^{\b}\neq y,z$,
or $y^{8^{-1}(\lambda^3-\lambda^2+7\lambda+1)}
z^{2^{-1}(1+\lambda)}$, contrary to Eq~(\ref{eq=cover2-eqS2}).

Again by Eq~(\ref{eq=cover2-eqS2}), we have $y^{\b}=y$ and $z^{\b}=(y^{\b}y^{-1})^{2^{-1}(\lambda-1)}z=z$.
Since $(xy)^{\b}=x^{\b}y\in S$, it is easy to check that
$(xy)^{\b}=xy$ and thus $x^{\b}=x$. Hence $\b$ is the identity automorphism of $H$
and $|\Aut(H,S)|=1$. By Observation, $|N_A(R(H))|=10|H|$.
\hfill$\blacksquare$

\section{Cyclic covers\label{s5}}

In this section, we classify connected symmetric cyclic covers of
connected pentavalent symmetric graphs of order twice a prime.
Denote by $K_{6,6}-6K_2$ the complete bipartite graph
of order 12 minus a one-factor and by $\mathbf{I}_{12}$ the Icosahedron graph.
Edge-transitive cyclic covers
of $K_6$ were classified in~\cite[Theorem~1.1]{PHXD}, and by~\cite[Line 20, pp.40]{PHXD}, such graphs have order 12 and thus isomorphic to $K_{6,6}-6K_2$ or $\mathbf{I}_{12}$ by~\cite[Proposition~2.7]{GHS} (note that the graph $\mathbf{I}_{12}$ is missed in~\cite[Theorem~1.1]{PHXD}).

\begin{theorem}\label{the=cycliccover}
Let $\G$ be a connected pentavalent symmetric
graph of order $2p$ for a prime $p$, and let $\widetilde{\G}$
be a connected symmetric $\mz_n$-cover of $\G$ with $n\geq2$.
Then $\widetilde{\G}\cong K_{6,6}-6K_2$, $\mathbf{I}_{12}$, $\mathcal{CD}_{np}$, or $\mathcal{CGD}_{mp^e\times p}^{i}$ for $1\leq i\leq 5$ with $n=mp^e$, $(m,p)=1$, $5\di (p-1)$ and $e\geq 1$, which are defined in Examples~\ref{exam=n}, \ref{exam=n*p}, \ref{exam=n*p*p-1} and \ref{exam=n*p*p-2}.
\end{theorem}

\f {\bf Proof:} By Proposition~\ref{prop=2p}, $\G\cong K_{6}$ for $p=3$,
$K_{5,5}$ for $p=5$, or $\mathcal{CD}_{p}$ for $5\ |\ (p-1)$. If $\G\cong K_6$ then
$\widetilde{\G}\cong K_{6,6}-6K_2$ or $\mathbf{I}_{12}$ by~\cite[Theorem~1.1]{PHXD} (also
see the proof in~\cite[Theorem~3.6]{AHK}). In the following, we assume that
$p\geq5$. Let $A=\Aut(\widetilde{\G})$.

Let $K=\mz_n$ and $F=N_{A}(K)$.
Since $\widetilde{\G}$ is a symmetric $K$-cover of $\G$,
$F$ is arc-transitive on $\widetilde{\G}$ and $F/K$ is
arc-transitive on $\widetilde{\G}_{K}=\G$.
Let $B/K$ be a minimal arc-transitive
subgroup of $F/K$.
By Proposition~\ref{prop=2p}, $B/K\cong D_{p}\rtimes\mz_5$ for $p>11$; by MAGMA~\cite{BCP}, $B/K\cong D_{11}\rtimes\mz_5$ for $p=11$, and $B/K\cong \mz_5^2\rtimes\mz_2$, $\mz_5^2\rtimes\mz_4$ or $\mz_5^2\rtimes\mz_8$ for $p=5$.
Each minimal normal subgroup of $B/K$ is isomorphic to
$\mz_p$ or $\mz_5^2$ with $p=5$ and $B/K\cong\mz_5^2\rtimes\mz_8$.
Clearly, $B$ is arc-transitive on $\widetilde{\G}$
and $B/K$ is non-abelian.

Set $C=C_B(K)$. Since $K$ is abelian, $K\leq Z(C)\leq C$, where $Z(C)$ is the center of $C$.
Suppose $K=C$. Then $B/K=B/C\lesssim \Aut(K)\cong\mz_n^*$, which forces that $B/K$ is abelian,
a contradiction. Hence $K<C$ and $1\neq C/K\unlhd B/K$. It follows that  $C/K$ contains a minimal normal subgroup of $B/K$, say $L/K$. Then $L\unlhd B$, $L\leq C\unlhd B$. Furthermore, $L/K\cong\mz_p$, or $L/K\cong\mz_5^2$ with $p=5$ and $B/K\cong\mz_5^2\rtimes\mz_8$.

Clearly, $L$ and $L/K$ have two orbits on $V(\widetilde{\G})$ and $V(\widetilde{\G}_{K})$, and $\widetilde{\G}$ and $\widetilde{\G}_{K}$ are bipartite graphs with the two orbits of  $L$ and $L/K$ as their bipartite sets, respectively. Since $K\leq Z(C)$ and $L\leq C$, $K\leq Z(L)$.

First, assume $L/K\cong\mz_p$. Since $K\leq Z(L)$, $L$ is abelian, and so $L\cong\mz_{np}$ or $\mz_n\times\mz_p$ with $p\ |\ n$.
For the latter, $L\cong\mz_m\times\mz_{p^e}\times\mz_p$
with $n=mp^e$, $(m,p)=1$ and $e\geq1$. Since $L/K$ is semiregular on $V(\widetilde{\G}_{K})$, $L$ is semiregular on $V(\widetilde{\G})$ and thus $\widetilde{\G}$ is
a bi-Cayley graph over $L$. Noting that $L\lhd B$, we have that
$N_{A}(L)$ is
arc-transitive on $\widetilde{\G}$, forcing that $\widetilde{\G}\cong\BiCay(L,\emptyset,\emptyset, S)$ for some subset $S\subseteq L$.
Recall that $p\geq5$. By Lemmas~\ref{lem=cycliccover}-\ref{lem=cover2},
$\widetilde{\G}\cong\mathcal{CD}_{np}$ or $\mathcal{CGD}_{mp^e\times p}^{i}$ ($1\leq i\leq 5$), as required.

Now, assume $L/K\cong\mz_5^2$. Then $p=5$ and $B/K\cong\mz_5^2\rtimes\mz_8$.
Since $K\leq Z(L)$ and $K=\mz_n$, $L=P\times H$, where $P$ and $H$ are the Sylow 5-subgroup
and the Hall $5'$-subgroup of $L$, respectively. Note that $H\leq K$ is abelian, but $P$ may not.
Since $(L/K)_{v^K}\cong\mz_5$, we have $L_v=P_v\cong\mz_5$, where $v\in V(\widetilde{\G})$ and $v^K$ is an orbit of $K$ on $V(\widetilde{\G})$ containing $v$.
Note that $P\unlhd B$ as $P$ is characteristic in
$L$ and $L\unlhd B$. By Proposition~\ref{prop=atlesst3orbits},
$P$ has at most two orbits on $V(\widetilde{\G})$ because $P_v\neq1$, and
since $L$ has exactly two orbits on $V(\widetilde{\G})$,
$P$ and $L$ have the same orbits. It follows that $L=PL_v=PP_v=P$, forcing
that $H=1$ and $K$ is a $5$-group.

Suppose $|K|=5^t$ with $t\geq 2$. Since $K$ is cyclic, $K$ has a characteristic subgroup $N$ such that $|K/N|=25$, and since $K\unlhd B$, $N\unlhd B$.
By Proposition~\ref{prop=atlesst3orbits},
$\widetilde{\G}_{N}$ is a connected pentavalent $B/N$-arc-transitive graph of order $10|K|/|N|=250$,
and by Example~\ref{exam=250}, $\G_N\cong\mathcal{CGD}_{5^3}$.
Since $B/K\cong\mz_5^2\rtimes\mz_8$ and $|K/N|=5^2$, all Sylow $2$-subgroups of $B/N$ are isomorphic to $\mz_8$ and $|B/N|=8\cdot 5^4$.
However, by MAGMA~\cite{BCP}, $\Aut(\mathcal{CGD}_{5^3})$ has no arc-transitive subgroup of order $8\cdot 5^4$ that has a Sylow $2$-subgroup isomorphic to $\mz_8$, a contradiction.

Since $K\not=1$, we have $|K|=5$ and $|V(\widetilde{\G})|=10|K|=50$.
By Example~\ref{exam=n*p*p-1}, $\widetilde{\G}\cong\mathcal{CGD}_{5\times 5}^{4}$, as required.  \hfill$\blacksquare$

\section{Dihedral covers}

In this section, we aim to classify symmetric dihedral covers
of connected pentavalent symmetric graphs of order twice a prime.
First, we introduce four graphs which are from~\cite{P}.

\begin{exam}\label{exam=G48} {\rm
Let
$\mathbf{I}_{12}^{(2)}=\Cay(D_{12},\{b,ba,ba^2,ba^{4},ba^{9}\})$
and $\mathcal{G}_{48}=\Cay(D_{24},\{b,ba,ba^3,ba^{11},ba^{20}\})$
be two Cayley graphs on the dihedral groups $D_{12}=\lg a,b~|~a^{12}=b^2=1,a^b=a^{-1}\rg$
and $D_{24}=\lg a,b~|~a^{24}=b^2=1,a^b=a^{-1}\rg$, respectively.
By MAGMA~\cite{BCP}, $\Aut(\mathbf{I}_{12}^{(2)})\cong A_5\rtimes D_4$ and
$\Aut(\mathcal{G}_{48})\cong {\rm SL}(2,5)\rtimes D_4$, and their vertex stabilizers are isomorphic
to $F_{20}$.

}

\end{exam}

\begin{exam}{\rm \label{exam=G60}
Let $\mathcal{G}_{60}=\Cay(A_5,\{(1\ 4)(2\ 5), (1\ 3)(2\ 5), (1\ 3)(2\ 4), (2\ 4)(3\ 5), (1\ 4)(3\ 5)\})$.
By MAGMA~\cite{BCP}, $\mathcal{G}_{60}$ is a connected pentavalent symmetric graph of order $60$ and $\Aut(\mathcal{G}_{60})\cong A_5\times D_5$ with vertex
stabilizer isomorphic to $D_5$.

}
\end{exam}

\begin{exam}{\rm \label{exam=G120}
Let $G$ be a subgroup of $S_7$ generated by $a=(1\ 4)(2\ 5)(6\ 7)$,
$b=(1\ 3)(2\ 5)(6\ 7)$, $c=(1\ 3)(2\ 4)(6\ 7)$, $d=(2\ 4)(3\ 5)(6\ 7)$
and $e=(1\ 4)(3\ 5)(6\ 7)$, and define
$\mathcal{G}_{120}=\Cay(G,\{a,b,c,d,e\})$.
By MAGMA~\cite{BCP}, $G\cong A_5\times\mz_2$ and $\mathcal{G}_{120}$ is a connected pentavalent symmetric graph of order $120$. Moreover, $\Aut(\mathcal{G}_{120})\cong A_5\times D_{10}$
with vertex stabilizer isomorphic to $D_5$.
}
\end{exam}

A list of all pentavalent $G$-arc-transitive graphs on up to 500 vertices with the vertex stabilizer $G_v\cong \mz_5$, $D_5$ or $F_{20}$ was given in magma code by Poto\v cnik~\cite{P}.
Based on this list, we have the following lemma.

\begin{lem}\label{lem=dihedral}
Let $\G$ be a $G$-arc-transitive graph of order $24$, $48$, $60$, $120$ or $240$
with vertex stabilizer $G_v\cong \mz_5$, $D_5$ or $F_{20}$ for some $G\leq\Aut(\G)$ and $v\in V(\G)$.
Then $\G$ is a connected symmetric dihedral cover of $K_6$ if and only if $\G\cong \mathbf{I}_{12}^{(2)}$, $\mathcal{G}_{48}$, $\mathcal{G}_{60}$ or $\mathcal{G}_{120}$.
\end{lem}

\f {\bf Proof:} To show the necessity, let $\G$ be a connected symmetric dihedral cover of $K_6$. Then $\Aut(\G)$ has an arc-transitive subgroup having a normal dihedral subgroup of order $|V(\G)|/6$. Since $\G$ is $G$-arc-transitive with $G_v\cong \mz_5$, $D_5$ or $F_{20}$, by~\cite{P} $\G$ is isomorphic to one of the seven graphs: three graphs of order 24, 48 and 60 respectively, two
graphs of order 120 and two graphs of order 240.
For the orders 24, 48 and 60, $\G\cong\mathbf{I}_{12}^{(2)}$, $\mathcal{G}_{48}$
 or $\mathcal{G}_{60}$ by Examples~\ref{exam=G48} and \ref{exam=G60}.
For the order 120, by MAGMA~\cite{BCP} one graph is isomorphic to $\mathcal{G}_{120}$
and the other has no arc-transitive group of automorphisms
having a normal dihedral subgroup of order $20$; in this case $\G\cong \mathcal{G}_{120}$. For the order 240, again by MAGMA~\cite{BCP} none of the two graphs
has an arc-transitive group of automorphisms
having a normal dihedral subgroup of order $40$.

Now, we show the sufficiency.
By MAGMA~\cite{BCP}, $\Aut(\mathbf{I}_{12}^{(2)})$ has a normal subgroup
$N\cong D_2$. Clearly, $N$ has more than two orbits on $V(\mathbf{I}_{12}^{(2)})$, and by Proposition~\ref{prop=atlesst3orbits},
the quotient graph $(\mathbf{I}_{12}^{(2)})_N$
is a connected pentavalent symmetric graph of order 6, that is, the complete graph $K_6$. Thus  $\mathbf{I}_{12}^{(2)}$ is a $D_2$-cover of $K_6$. Similarly, one may show that $\mathcal{G}_{48}$, $\mathcal{G}_{60}$ or $\mathcal{G}_{120}$ is a symmetric $D_3$-, $D_5$- or $D_{10}$-cover of
$K_6$, respectively.
\hfill$\blacksquare$

Now, we are ready to classify symmetric dihedral covers
of connected pentavalent symmetric graphs of order $2p$ for any prime $p$.
Clearly, we have $p\geq 3$.

\begin{theorem}\label{the=Dihedral-cover}
Let $\G$ be a connected pentavalent symmetric graph of order $2p$ with $p$ a prime,
and let $\widetilde{\G}$ be a connected symmetric $D_n$-cover of $\G$ with $n\geq2$.
Then $\widetilde{\G}\cong
\mathbf{I}_{12}^{(2)}$, $\mathcal{G}_{48}$,  $\mathcal{G}_{60}$ or
$\mathcal{G}_{120}$.
\end{theorem}

\f {\bf Proof:} Let $K=D_n$ and let $F$ be the fibre-preserving group.
Since $\widetilde{\G}$ is a symmetric $K$-cover of $\G$,
$F$ is arc-transitive on $\widetilde{\G}$ and $F/K$ is arc-transitive on $\widetilde{\G}_K=\G$.

Assume $n=2$. Then $|V(\widetilde{\G})|=2n\cdot |V(\G)|=8p$.
Recall that $p\geq 3$. By~\cite[Proposition~2.9]{GHS}, $\widetilde{\G}\cong
\mathbf{I}_{12}^{(2)}$ or a graph $\mathcal{G}_{248}$ of order 248 with
$\Aut(\mathcal{G}_{248})=\PSL(2,31)$.
Since $\PSL(2,31)$ has no proper subgroup of order divisible by $248$ by MAGMA~\cite{BCP},
$\Aut(\widetilde{\G})$ is the unique arc-transitive group of automorphisms of $\widetilde{\G}$, that is, $F\cong \PSL(2,31)$. It implies that $\widetilde{\G}\not\cong \mathcal{G}_{248}$ because $F$ has no normal subgroup isomorphic to $D_n$. Hence $\widetilde{\G}\cong\mathbf{I}_{12}^{(2)}$.

Assume $n>2$. Let $\mz_n$ be the cyclic subgroup of $K=D_{n}$
of order $n$. Then ${\mz_n}$ is characteristic in $K$ and so ${\mz_n}\unlhd F$ as $K\unlhd F$.
By Propositions~\ref{prop=atlesst3orbits},
$\widetilde{\G}_{\mz_n}$ is a connected pentavalent $F/{\mz_n}$-arc-transitive graph of order $4p$, and by~\cite[Proposition~2.7]{GHS},
$\widetilde{\G}_{\mz_n}\cong \mathbf{I}_{12}$ or $K_{6,6}-6K_2$. Thus $\widetilde{\G}$ is a symmetric ${\mz_n}$-cover of $K_{6,6}-6K_2$ or $\mathbf{I}_{12}$.
Note that $|V(\widetilde{\G})|=12n$.

Let $\widetilde{\G}_{\mz_n}\cong K_{6,6}-6K_2$. Since each
minimal arc-transitive subgroup of $\Aut(K_{6,6}-6K_2)$
is isomorphic to $A_5\times \mz_2$ or $S_5$ by MAGMA~\cite{BCP},
$F/{\mz_n}$ has an arc-transitive subgroup
$B/{\mz_n}=A_5\times \mz_2$ or $S_5$. It follows that $|B_v|=10$ for $v\in V(\widetilde{\G})$, and form Proposition~\ref{prop=stabilizer} that $B_v\cong D_5$.
In particular, $B$ is arc-transitive on $\widetilde{\G}$ and $B/{\mz_n}$ has a normal
subgroup $M/{\mz_n}=A_5$, which is edge-transitive on $\widetilde{\G}_{\mz_n}$ and has exactly two orbits on $V(\widetilde{\G}_{\mz_n})$.
Thus $M\unlhd B$ is edge-transitive and has two orbits on $V(\widetilde{\G})$. Since $|B:M|=2$, we have $M_v\cong D_5$.

Clearly, $\mz_n\leq C_M(\mz_n)$. If $\mz_n=C_M(\mz_n)$, then
$A_5=M/\mz_n=M/C_M(\mz_n)\leq\Aut(\mz_n)=\mz_n^*$, which is impossible.
Hence $\mz_n$ is a proper subgroup of $C_M(\mz_n)$, and since
$\Mult(A_5)=\mz_2$, Lemma~\ref{lem=Mult} implies that
either $M=M'\times {\mz_n}=A_5\times \mz_n$ or $M=M'{\mz_n}={\rm SL}(2,5)\mz_n$ with $M'\cap {\mz_n}\cong\mz_2$. In particular, $M/M'$ is cyclic.
Since $M'$ is characteristic in $M$ and
$M\unlhd B$, we have $M'\unlhd B$.
If $M'$ has at least three orbits on $V(\widetilde{\G})$, by Proposition~\ref{prop=atlesst3orbits}, $M'$ is semiregular on $V(\widetilde{\G})$ and  $\widetilde{\G}_{M'}$ is a connected pentavalent $B/M'$-arc-transitive graph.
The stabilizer of $\a\in V(\widetilde{\G}_{M'})$ in $M/M'$ is isomorphic to $M_v\cong D_5$, but this is impossible because $M/M'$ is cyclic. Thus $M'$ has at most two orbits on $V(\widetilde{\G})$ and so $|V(\widetilde{\G})|\di 2|M'|$, that is, $6n\di |M'|$.
If $M=A_5\times \mz_n$, then $M'=A_5$ and $6n\di |M'|$ implies that $n=5$ or 10 as $n>2$. It follows that $|V(\widetilde{\G})|=60$ or 120. Since $B_v\cong D_5$, we have
$\tilde{\G}\cong \mathcal{G}_{60}$ or $\mathcal{G}_{120}$ by Lemma~\ref{lem=dihedral}.
If $M={\rm SL}(2,5)\mz_n$ with $M'={\rm SL}(2,5)$ and
${\rm SL}(2,5)\cap \mz_n\cong\mz_2$, then
$n$ is even and $6n\di |M'|$ implies that $n=4$, 10 or 20. It follows that $|V(\widetilde{\G})|=48$, 120 or 240, and from Lemma~\ref{lem=dihedral} that $\tilde{\G}\cong \mathcal{G}_{48}$ or $\mathcal{G}_{120}$.

Let $\widetilde{\G}_{\mz_n}\cong \mathbf{I}_{12}$.
By MAGMA~\cite{BCP}, under conjugation
$\Aut(\mathbf{I}_{12})$ has only one minimal arc-transitive subgroup isomorphic to $A_5$, and so $F/{\mz_n}$ has an arc-transitive subgroup
$B/{\mz_n}\cong A_5$. By a similar argument as the previous paragraph,
$B=B'{\mz_n}$ and $B'\cap {\mz_n}\lesssim \Mult(A_5)$ by Lemma~\ref{lem=Mult},
forcing that either $B=B'\times {\mz_n}=A_5\times \mz_n$ or $B=B'{\mz_n}={\rm SL}(2,5)\mz_n$ with ${\rm SL}(2,5)\cap \mz_n\cong\mz_2$. Furthermore, $B$ is arc-transitive on $\widetilde{\G}$ with $B_v\cong \mz_5$ for $v\in V(\widetilde{\G})$, and $B/B'$ is cyclic. If $B'$ has more than two orbits on $V(\G)$, then $\widetilde{\G}_{B'}$ is a connected pentavalent $B/B'$-arc-transitive graph by Proposition~\ref{prop=atlesst3orbits}, which is impossible because $B/B'$ is abelian. Thus $B'$ has at most two orbits on $V(\widetilde{\G})$ and so $12n\di 2|B'|$.
If $B=A_5\times \mz_n$, then $B'\cong A_5$, and $12n\di 2|B'|$ implies that $n=5$ or $10$. It follows that $|V(\widetilde{\G})|=60$ or $120$. Since $B_v\cong\mz_5$, we have $\widetilde{\G}\cong\mathcal{G}_{60}$
or $\mathcal{G}_{120}$ by Lemma~\ref{lem=dihedral}. If $B={\rm SL}(2,5)\mz_n$ with ${\rm SL}(2,5)\cap \mz_n\cong\mz_2$, then $B'\cong{\rm SL}(2,5)$ and $n$ is even. Since $12n\di 2|B'|$, we have  $n=4$, 10 or 20, and so $|V(\widetilde{\G})|=48$, 120 or 240. It follows from Lemma~\ref{lem=dihedral} that
$\widetilde{\G}\cong \mathcal{G}_{48}$
or $\mathcal{G}_{120}$. \hfill$\blacksquare$

\section{Full automorphism groups of covers}

Let $\G$ be a symmetric $D_n$- or $\mz_n$-cover of
a connected symmetric pentavalent graph of order $2p$,
where $n\geq2$ is an integer and $p$ is a prime.
In this section, we aim to determine the
full automorphism group of $\G$.
For $D_n$, by Theorem~\ref{the=Dihedral-cover},
$\G\cong\mathbf{I}_{12}^{(2)}$, $\mathcal{G}_{48}$,  $\mathcal{G}_{60}$
or $\mathcal{G}_{120}$ and by Examples~\ref{exam=G48}-\ref{exam=G120},
$\Aut(\G)$ is known.
For $\mz_n$, by Theorem~\ref{the=cycliccover},
$\G\cong K_{6,6}-6K_2$, $\mathbf{I}_{12}$, $\mathcal{CD}_{np}$ (see Example~\ref{exam=n}), or $\mathcal{CGD}_{mp^e\times p}^{i}$ with $1\leq i\leq 5$ (see Examples~\ref{exam=n*p}, \ref{exam=n*p*p-1} and \ref{exam=n*p*p-2}). In particular, for the graph $\mathcal{CGD}_{mp^e\times p}^{i}$, we have $mp^e=n$ and $m$ is given by

\parbox{8cm}{
\begin{eqnarray*}
&& m=5^tp_1^{e_1}p_2^{e_2}\cdots p_s^{e_s} \mbox{\ \ \ s.t. \ \ \ $t\leq 1$, $s\geq 0$, $e_j\geq1$, $5\ |\ (p_j-1)$ for $0\leq j\leq s$,}
\end{eqnarray*}}\hfill
\parbox{1cm}{
\begin{eqnarray}
\label{eq=m}
\end{eqnarray}}

\f where $m,p,e$ satisfy the conditions
as listed in the second collum in Table~\ref{table=3}.
Note that $m$ is odd by Eq~(\ref{eq=m}). By MAGMA~\cite{BCP},
$\Aut(K_{6,6}-6K_2)=S_6\times\mz_2$ and
$\Aut(\mathbf{I}_{12})=A_5\times\mz_2$, and by Example~\ref{exam=n},
$\Aut(\mathcal{CD}_{np})=D_{np}\rtimes\mz_5$.
Hence we only need to determine the full automorphism
groups of $\mathcal{CGD}_{mp^e\times p}^{i}$ for $1\leq i\leq 5$. All theses graphs are connected symmetric cyclic covers of pentavalent symmetric graphs of order $2p$ except $\mathcal{CGD}_{mp\times p}^{4}$ with $5\di (p+1)$, which are connected symmetric bi-Cayley graphs over $\mz_{mp}\times\mz_p$.

\begin{theorem}\label{lem=automorphism}
$\Aut(\mathcal{CGD}_{mp^e\times p}^{i})$ for $1\leq i\leq 5$
is isomorphic to one group listed in Table~\ref{table=3}.

\begin{table}[ht]

\begin{center}
\begin{tabular}{|l|l|l|}

\hline
$\G$              & Conditions: $(m,p)=1$, $m$: Eq~(\ref{eq=m})              & $\Aut(\G)$   \\
\hline
$\mathcal{CGD}_{mp^e\times p}^{i}~(i=1,2,3)$ & $5\ |\ (p-1)$ and $e\geq2$ & $\Dih(\mz_{mp^e}\times\mz_p)\rtimes\mz_5$ \\
\hline
\multirow{3}{*}{$\mathcal{CGD}_{mp\times p}^{4}$}  & $m\neq 1,5$, and $p=5$ or $5\ |\ (p\pm1)$ &$\Dih(\mz_{mp}\times\mz_p)\rtimes \mz_5$\\
\cline{2-3}
& $m=1$ or $5$, and $5\ |\ (p\pm1)$
& $\Dih(\mz_{mp}\times\mz_p)\rtimes D_5$ \\
\cline{2-3}
 &$m=1$ and $p=5$ & $(\Dih(\mz_5^2)\times F_{20}).\mz_4$\\
\hline
$\mathcal{CGD}_{mp\times p}^{5}$  & $5\ |\ (p-1)$ & $\Dih(\mz_{mp}\times\mz_p)\rtimes \mz_5$\\
\hline
\end{tabular}
\end{center}
\vskip -0.5cm
\caption{{\small Full automorphism groups of $\mathcal{CGD}_{mp^e\times p}^{i}$
for $1\leq i\leq5$}}\label{table=3}
\end{table}

\end{theorem}

\f {\bf Proof:} Let $\G=\mathcal{CGD}_{mp^e\times p}^{i}$ for $1\leq i\leq 5$
and $A=\Aut(\G)$. For $(m,p)=(1,5)$, we have $\G=\mathcal{CGD}_{5\times 5}^{4}$ and by~\cite[Theorem~4.3~(1)]{FZL}, $\Aut(\G)\cong(\Dih(\mz_5^2)\rtimes F_{20}).\mz_4$. In what follows we assume that $(m,p)\neq (1,5)$.
By Examples~\ref{exam=n*p}, \ref{exam=n*p*p-1} and \ref{exam=n*p*p-2}, $A$ has an arc-transitive subgroup $F$ isomorphic to
$\Dih(\mz_{mp^e}\times\mz_p)\rtimes\mz_5$ for $\mathcal{CGD}_{mp^e\times p}^{i}$ ($i=1,2,3$),  $\Dih(\mz_{mp}\times\mz_p)\rtimes \mz_5$ for
$\mathcal{CGD}_{mp\times p}^{4}$ with $m\neq 1,5$ and $p=5$ or $5\di (p\pm 1)$,
$\Dih(\mz_{mp}\times\mz_p)\rtimes D_5$ for $\mathcal{CGD}_{mp\times p}^{4}$ with $m=1$ or $5$ and $5\di (p\pm 1)$, and $\Dih(\mz_{mp}\times\mz_p)\rtimes \mz_5$ for $\mathcal{CGD}_{mp\times p}^{5}$ with $5\di (p-1)$. Note that $F_v=\mz_5$ or $D_5$ for $v\in V(\G)$. Furthermore, $F$ has a normal semiregular subgroup $K=\mz_{mp^e}\times\mz_p$ having two orbits on $V(\G)$, and hence $\G$ is an F-arc-transitive bi-Cayley graph over $K$. By Lemmas~\ref{lem=cover1} and \ref{lem=cover2}, $|N_A(K)|=|F|$, implying that $N_A(K)=F$.
Note that $|F|=10|K|$ or $20|K|$, that is,  $|F|=10mp^{e+1}$ or $20mp^{e+1}$ with $p=5$ or $5\di (p\pm 1)$, and by Eq~(\ref{eq=m}), both $m$ and $|K|$ are odd. In particular, $|V(\G)|=2|K|=2mp^{e+1}$ is twice an odd integer.

Clearly, $K=\mz_{mp^e}\times\mz_p$ has a characteristic Hall $5'$-subgroup, say $H$. Then $H\unlhd F$ as $K\unlhd F$. If $H\neq K$, then $5\ |\ mp^{e+1}$ and $H$ has at least three orbits. For $p\neq5$, we have $5\ |\ m$, and since $5^2\nmid m$ by Eq~(\ref{eq=m}), we have
$|K:H|=5$. For $p=5$,
by Table~\ref{table=3}, $\G=\mathcal{CGD}_{mp\times p}^4$ with $(m,5)=1$ and $K=\mz_m\times\mz_p\times\mz_p$,
implying that $|K:H|=5^2$.
By Proposition~\ref{prop=atlesst3orbits},
$\G_H$ is a connected pentavalent $F/H$-arc-transitive
graph of order $2\cdot 5$ or $2\cdot 5^2$.
By Proposition~\ref{prop=2p} and Example~\ref{exam=n*p*p-1}, $\G_H\cong K_{5,5}$
or $\G_H\cong\mathcal{CGD}_{5\times 5}^{4}$. Since $|F|=10|K|$ or $20|K|$ and $|K|$ is odd, $H$ is the characteristic Hall $\{2,5\}'$-subgroup of $F$.
Thus we have the following claim.

\vskip 0.2cm
\f {\bf Claim:} $H$ is the characteristic Hall $\{2,5\}'$-subgroup of $F$, and we have $H=K$, or $|K:H|=5$ and $\G_H\cong K_{5,5}$, or $|K:H|=25$ and $\G_H\cong \mathcal{CGD}_{5\times 5}^{4}$.

\vskip 0.2cm

To finish the proof, we only need to show that $A=F$. Suppose to the contrary that $A\neq F$. Then $A$ has a subgroup $M$ such that $F$ is a maximal subgroup of $M$.
Since $F$ is arc-transitive on $\G$,
$M$ is arc-transitive, and since $N_A(K)=F$, we have  $K\ntrianglelefteq M$.

By the definitions of the graphs $\mathcal{CGD}_{mp^e\times p}^i$ ($1\leq i\leq5$) in Examples~\ref{exam=n*p}, \ref{exam=n*p*p-1} and \ref{exam=n*p*p-2}, $\G$ has the 6-cycle
$(1, h, a^{-r-1}b^{-\lambda-1}c^{-1},
ha^{-r}b^{-\lambda}c^{-1}, a^{-r}b^{-\lambda}c^{-1}, hab, 1)$ for $1\leq i\leq3$, and the $6$-cycle $(1,h,a^{-r-1}c^{-1},ha^{-r}bc^{-1},a^{-r}bc^{-1},hab,1)$
for $4\leq i \leq5$.
Suppose that $\G$ is $(M,4)$-arc-transitive. Then each 4-arc lies in a 6-cycle
in $\G$ and so $\G$ has diameter at most three. It follows that $|V(\G)|=2mp^{e+1}\leq 1+5+5\cdot 4+5\cdot4\cdot4=106$, that is, $mp^{e+1}\leq 53$.
Since $p=5$ or $5\di (p\pm 1)$ and $e+1\geq2$
(see the second collum of Table~\ref{table=3}),
we have $p=5$ and $m\leq2$. Since $m$ is odd, $(m,p)=(1,5)$, contrary to assumption. Thus $\G$ is at most $3$-arc-transitive, and by Proposition~\ref{prop=stabilizer},
we have $|M_v|\in\{5,10,20,40,60,80,120,720,1440,2880\}$.

Note that $|M:F|=|M_v:F_v|\in\{2, 4, 6, 8, 12, 16, 24, 72, 144, 288, 576\}$ because $M\not=F$ and $|F_v|=5$ or $10$.
Let $[M:F]$ be the set of right cosets of $F$ in $M$. Consider the action of $M$ on  $[M:F]$ by right multiplication, and let $F_M$ be the kernel of this action, that is, the largest normal subgroup of $M$ contained in $F$. Then $M/F_M$ is a primitive permutation group on $[M:F]$ because $F/F_M$ is maximal in $M/F_M$, and $(M/F_M)_F=F/F_M$, the stabilizer of $F\in [M:F]$ in $M/F_M$. It follows that $|M/F_M|=|M:F||F/F_M|$ and so  $|F/F_M|=|M/F_M|/|M:F|$. Since $|M:F|\in\{2, 4, 6, 8, 12, 16, 24, 72, 144, 288, 576\}$, by Lemma~\ref{lem=primitive permutation group} we have $M/F_M\leq {\rm AGL}(t,2)$ with $|M:F|=2^t$ and $1\leq t\leq 4$, or
$\soc(M/F_M)\cong \PSL(2,q)$, $\PSL(3,3)$ or $\PSL(2,r)\times \PSL(2,r)$ with
$|M:F|=q+1$, $144$ or $(r+1)^2$ respectively, where
$q\in\{5,7,11,23,71\}$ and $r\in \{11, 23\}$.

Suppose $M/F_M\leq {\rm AGL}(2,2)$ and $|M:F|=4$. Since a $2$-group cannot be primitive on $[M:F]$, we have $3\di |M/F_M|$ and so $3\di |M/F_M|/|M:F|=|F/F_M|$. Since $|F|=10mp^{e+1}$ or $20mp^{e+1}$ with $p=5$ or $5\di (p\pm 1)$, we have $3\di m$, which is impossible by Eq~(\ref{eq=m}). Thus $M/F_M\not\leq {\rm AGL}(2,2)$. Similarly, since $7\nmid m$, we have $M/F_M\not\leq {\rm AGL}(3,2)$, and if $M/F_M\leq {\rm AGL}(4,2)$, then $M/F_M$ is a $\{2,5\}$-group. Furthermore, $\soc(M/F_M)\not\cong \PSL(2,q)$, $\PSL(3,3)$ or $\PSL(2,r)\times \PSL(2,r)$ for $q\in\{7,23,71\}$ and $r=23$ because otherwise one of $7,23,13,23$ is a divisor of $m$. It follows that
$M/F_M\cong \mz_2$ with $|M:F|=2$, $M/F_M\leq {\rm AGL}(4,2)$ with $|M:F|=2^4$ and $M/F_M$ a $\{2,5\}$-group, $\soc(M/F_M)\cong \PSL(2,q)$ with $|M:F|=q+1$ and $q\in\{5,11\}$, or $\soc(M/F_M)\cong\PSL(2,11)\times \PSL(2,11)$ with $|M:F|=144$.

First assume that $M/F_M\cong\mz_2$ with $|M:F|=2$. Then $F\unlhd M$ and $H\unlhd M$ as $H$ is characteristic in $F$ by Claim. Let $C=C_M(H)$. Since $K$ is abelian, $H\leq K\leq C$.  Let $P$ be a Sylow 5-subgroup of $C$ containing the unique Sylow $5$-subgroup of $K$. Since $H$ is the Hall $5'$-group of $K$, $K\leq HP=H\times P$.
Clearly, $HP/H$ is a Sylow 5-subgroup of $C/H$.
Recall that $|F/K|\di 20$ and $|K/H|\di 25$ (see Claim).
Since $|M|=2|F|$, we have $|M/H|\di 2^3\cdot 5^3$, and by Sylow theorem, $M/H$ has a normal Sylow $5$-subgroup. In particular, $C/H$ has a normal Sylow $5$-subgroup, that is, $HP/H\unlhd C/H$. This implies $H\times P\unlhd C$, and since $C\unlhd M$ and $P$ is characteristic in $C$, we have $P\unlhd M$.
Since $(m,p)\neq (1,5)$ and $|V(\G)|=2mp^{e+1}$,
$P$ has at least three orbits on $V(\G)$. By Proposition~\ref{prop=atlesst3orbits}, $P$ is semiregular on $V(\G)$. Thus $|P|\di |V(\G)|$ and $|P|\di |K|$. It follows that $|HP|=|H||P|\di |K|$, and since $K\leq HP$, we have $K=HP\unlhd M$, a contradiction.

Assume that $M/F_M\leq {\rm AGL}(4,2)$ with $|M:F|=2^4$ and $M/F_M$ a $\{2,5\}$-group. Then $M/F_M$ has a regular normal subgroup of order $2^4$, say $L/F_M$, and hence $L\unlhd M$, $2^4\di |L|$ and $5\di |M:L|$.
If $L$ is semiregular then $2^4\di |V(\G)|=2mp^{e+1}$, which is impossible. Thus $L$ is not semiregular, and so $5\di |L_v|$. By Proposition~\ref{prop=atlesst3orbits}, $L$ has one or two orbits, yielding that $|L|=|V(\G)||L_v|$ or $|L|=|V(\G)||L_v|/2$. Since $|M|=|V(\G)||M_v|$, we have $|M:L|=|M_v:L_v|$ or $2|M_v:L_v|$, and since $5^2\nmid |M_v|$, we have $5\nmid |M:L|$, a contradiction.

Assume that $\soc(M/F_M)\cong \PSL(2,5)$ with $|M:F|=6$.
Then $M/F_M=\PSL(2,5)$ or $\PGL(2,5)$, and $|F/F_M|=|M/F_M|/|M:F|=10$ or $20$.
Since $H$ is the unique normal Hall $\{2,5\}'$-subgroup of $F$,
we have $H\leq F_M$ and so $H$ is characteristic in $F_M$. This implies $H\unlhd M$ because $F_M\unlhd M$.
Since $M/F_M\cong (M/H)/(F_M/H)$, $M/H$ is insolvable, and since $K\ntrianglelefteq M$, we have $H\neq K$. By Claim, $\G_H\cong K_{5,5}$ or $\mathcal{CGD}_{5\times5}^{4}$.
If $\G_H\cong \mathcal{CGD}_{5\times5}^{4}$ then $\Aut(\G_H)\cong (\Dih(\mz_5^2)\times F_{20}).\mz_4$ is solvable and so $M/H$ is solvable, a contradiction. If  $\G_H\cong K_{5,5}$ then as $\Aut(K_{5,5})=(S_5\times S_5)\rtimes\mz_2$, it is easy to show that each insolvable arc-transitive group of $\Aut(K_{5,5})$ contains $A_5\times A_5$ (this is also easily checked by MAGMA~\cite{BCP}), and so $|M/H|\geq 2\cdot 60^2$.
Noting that $F_M$ is semiregular on $V(\G)$,  we have $|F_M|\di |K|$. By Claim, $|K:H|\di 5^2$, and hence $|F_M:H|\di 5^2$. It follows that $|M/F_M|=|M/H|/|F_M/H|\geq 2\cdot 60^2/5^2>|\PGL(2,5)|$, a contradiction.

Assume that $L/F_M:=\soc(M/F_M)\cong\PSL(2,11)$ with $|M:F|=12$.
Then $M/F_M=\PSL(2,11)$ or $\PGL(2,11)$, and $|F/F_M|=|M/F_M|/|M:F|=55$ or 110. Moreover, $L\unlhd M$
and $K\leq L$ as $|K|$ is odd and $|M:L|\leq 2$.
Since $11\di |L/F_M|$, $F_M$ has at least three orbits on $V(\G)$, and by Proposition~\ref{prop=atlesst3orbits} $F_M$ is semiregular and $\G_{F_M}$ is a pentavalent $F/F_M$-arc-transitive graph. Thus $|F_M|\di |V(\G)|$ and $|V(\G_{F_M})|$ is even. Since $|V(\G_{F_M})|=|V(\G)|/|F_M|=2|K|/|F_M|$, $|F_M|$ is odd and  $|F_M|\di |K|$.

Set $N=H\cap F_M$.
Since $H$ is the characteristic Hall $\{2,5\}'$-subgroup of $F$, $N$ is the characteristic Hall $5'$-subgroup of $F_M$, and thus $N\unlhd M$ as $F_M\unlhd M$. Hence
$F_M/N$ is a 5-subgroup. By Claim, $5^3\nmid |K|$, and since $|F_M|\di |K|$, we have $5^3\nmid |F_M|$, that is, $|F_M/N|\di 25$.
Thus $F_M/N$ is abelian, and $\Aut(F_M/N)$ is cyclic or $\Aut(F_M/N)\cong{\rm GL}(2,5)$. If $F_M/N=C_{L/N}(F_M/N)$, then $\PSL(2,11)\cong L/F_M\cong (L/N)/(F_M/N)\lesssim \Aut(F_M/N)$, which is clearly impossible. Thus $F_M/N$ is a proper subgroup of $C_{L/N}(F_M/N)$, and since $\Mult(\PSL(2,11))\cong\mz_2$, Lemma~\ref{lem=Mult} implies that $L/N=(L/N)'\times F_M/N$ with $(L/N)'\cong \PSL(2,11)$. Since $|V(\G_N)|=|V(\G)|/|N|=2|K|/|N|$ with $|K|$ odd, $(L/N)'\cong\PSL(2,11)$ cannot be semiregular on $V(\G_N)$, implying that $5\di |(L/N)'_\a|$ for $\a\in V(\G_N)$. By Proposition~\ref{prop=atlesst3orbits}, $(L/N)'$ has at most two orbits on $V(\G_N)$, and so $|(L/N)|/|(L/N)'|=|V(\G_N)||(L/N)_\a|/(|V(\G_N)||(L/N)'_\a|)=
|(L/N)_\a|/|(L/N)'_\a|$ or $2|(L/N)_\a|/|(L/N)'_\a|$ is a $\{2,3\}$-group. Since $|(L/N)/(L/N)'|=|F_M/N|$ is a $5$-group, we have $|F_M/N|=1$, that is, $L/N=(L/N)'\cong \PSL(2,11)$.

Since $K\ntrianglelefteq M$ and $N\unlhd M$, we have $N\not=K$, and since $K\leq C_L(N)$ and $|N|$ is odd, Lemma~\ref{lem=Mult} implies $L=L'\times N$ with $L'\cong\PSL(2,11)$. Note that $L'\unlhd M$. Since $\G$ has order twice an odd integer, $L'$ cannot be semiregular on $\G$, yielding $5\di |L'_v|$. By Proposition~\ref{prop=atlesst3orbits}, $L'$ has at most two orbits, and so $|\PSL(2,11)|=|L'|=|V(\G)||L'_v|$ or $|V(\G)||L'_v|/2$. Since $3\nmid |V(\G)|=2mp^{e+1}$, we have $|V(\G)|=22$, contrary to the fact that $e+1\geq 2$.

Assume that $L/F_M:=\soc(M/F_M)\cong\PSL(2,11)\times\PSL(2,11)$
with $|M:F|=144$. Then there exists $L_1/F_M\unlhd L/F_M$ such that $L_1/F_M\cong\PSL(2,11)$ and $11\di |L:L_1|$.
Since $11\di |L:F_M|$, $F_M$ has at least three orbits and so $\G_{F_M}$ has order twice an odd integer. This implies that $L/F_M$ cannot be semiregular, and by  Proposition~\ref{prop=atlesst3orbits}, $L/F_M$ has one or two orbits. If $L/F_M$ has one orbit then $L_1/F_M$ is semiregular on $\G_{F_M}$ as  $11\di |L:L_1|$ implies that $L_1/F_M$ has at least three orbits, and so $4\di |V(\G_{F_M})|$, a contradiction. If $L/F_M$ has two  orbits then $\G_{F_M}$ is bipartite and $L/F_M$ is edge-transitive on $\G_{F_M}$. Furthermore,
$L_1/F_M$ fixes the bipartite sets setwise. Since $11\di |L:L_1|$, $L_1/F_M$ has at least two orbits on each bipartite set, and by~\cite[Proposition~2.4]{GHS}, $L_1/F_M$ is semiregular on $\G_{F_M}$. Since $L_1/F_M\cong \PSL(2,11)$, again we have the contradiction that $4\di |V(\G_{F_M})|$.
\hfill$\blacksquare$

\medskip

\f {\bf Acknowledgement:} This work was supported by the National Natural Science Foundation of China (11571035, 11231008, 11271012, 11371052) and by the 111 Project of China (B16002).

\end{document}